\documentclass[11pt]{article}
\usepackage{amscd}
\usepackage{amsmath}
\usepackage{latexsym}
\usepackage{amsfonts}
\usepackage{amssymb}
\usepackage{amsthm}
\usepackage{graphicx}

\theoremstyle{plain}
\newtheorem{theorem}{Theorem}
\newtheorem{proposition}[theorem]{Proposition}
\newtheorem{lemma}[theorem]{Lemma}
\newtheorem{corollary}[theorem]{Corollary}
\newtheorem{conjecture}[theorem]{Conjecture}
\theoremstyle{definition}
\newtheorem{definition}[theorem]{Definition}
\newtheorem{remark}[theorem]{Remark}

\numberwithin{equation}{section} \numberwithin{theorem}{section}

\begin{document}

 \title{ Global wellposedness and scattering for the focusing energy-critical nonlinear
Schr\"odinger equations of fourth order in the radial case}
\author{{Changxing Miao$^{\dag}$\ \ Guixiang Xu$^{\dag}$ \ \ and \ Lifeng Zhao $^{\ddag}$}\\
         {\small $^{\dag}$Institute of Applied Physics and Computational Mathematics}\\
         {\small P. O. Box 8009,\ Beijing,\ China,\ 100088}\\
         {\small $^\ddag$ Department of Mathematics, University of Science and Technology of China}\\
         {\small Hefei,\ China,\ 230026}\\
         {\small (miao\_changxing@iapcm.ac.cn, \ xu\_guixiang@iapcm.ac.cn, zhaolifengustc@yahoo.cn ) }\\
         \date{}
        }
        \maketitle

\begin{abstract} We consider the focusing energy-critical nonlinear
Schr\"odinger equation of fourth order $iu_t+\Delta^2
u=|u|^\frac{8}{d-4}u$. We prove that if a maximal-lifespan radial
solution $u: I\times\Bbb R^d\rightarrow\mathbb{C}$ obeys
$\displaystyle\sup_{t\in I}\|\Delta u(t)\|_{2}<\|\Delta W\|_{2}$,
then it is global and scatters both forward and backward in time.
Here $W$ denotes the ground state, which is a stationary solution of
the equation. In particular, if a solution has both energy and
kinetic energy less than those of the ground state $W$ at some point
in time, then the solution is global and scatters.
\end{abstract}
 \begin{center}
  \begin{minipage}{120mm}
 { \small {\bf Key Words:}
      {Focusing, Energy-critical, Radial, Fourth order Schr\"odinger equations, Global wellposedness, Scattering.}
   }\\
    { \small {\bf AMS Classification:}
      { 35Q40, 35Q55, 47J35.}
      }
      \end{minipage}
 \end{center}
\section{Introduction}
Fourth-order Schr\"odinger equations have been introduced by Karpman
\cite{Ka} and Karpman, Shagalov \cite{Ka2} to take into account the
role of small fourth-order dispersion terms in the propagation of
intense laser beams in a bulk medium with Kerr nonlinearity. Such
fourth-order Schr\"odinger equations are written as
\begin{equation}\label{gfnls}
i\partial_tu+\Delta^2u+\varepsilon\Delta u+f(|u|^2)u=0,
\end{equation}
where $\varepsilon=\pm1$ or $\varepsilon=0$, and $u: I\times\Bbb
R^d\rightarrow\Bbb C$ is a complex valued function. In this paper,
we will investigate the focusing energy-critical case when
$\varepsilon=0$, namely,
\begin{equation}\label{fnls} \left\{ \aligned
    iu_t +  \Delta^2 u  & = |u|^\frac{8}{d-4}u, \quad  \text{in}\  \mathbb{R}^d \times \mathbb{R}, d\geq5\\
     u(0)&=u_0(x), \quad \text{in} \ \mathbb{R}^d.
\endaligned
\right.
\end{equation} The name `energy-critical' refers to the fact that the scaling
symmetry
\begin{equation*} u(t,x)\mapsto
u_\lambda(t,x):=\lambda^\frac{d-4}{2}u(\lambda^4t, \lambda x)
\end{equation*}
leaves both the equation and the energy invariant. The energy of a
solution is defined by
\begin{equation*}
E(u(t))=\frac{1}{2}\int_{\Bbb R^d}|\Delta
u(t,x)|^2dx-\frac{d-4}{2d}\int_{\Bbb R^d}|u(t,x)|^\frac{2d}{d-4}dx
\end{equation*}
and is conserved under the flow. We refer to the Laplacian term in
the formula above as the kinetic energy and to the second term as
the potential energy.
\begin{definition}[Solutions.] Let $d\geq5$. A function $u: I\times\Bbb
R^d\rightarrow\Bbb C$ on a non-empty time interval $I\subset\Bbb R$
is a solution to (\ref{fnls}) if it lies in the class
$C_t^0\dot{H}_x^2(K\times\Bbb R^d)\cap
L_{t,x}^\frac{2(d+4)}{d-4}(K\times\Bbb R^d)$ for all compact
$K\subset I$, and obeys the Duhamel formula
\begin{equation*}
u(t)=e^{i(t-t_0)\Delta^2}u(t_0)-i\int_{t_0}^t
e^{i(t-\tau)\Delta^2}F(u(\tau))d\tau
\end{equation*}
for all $t,t_0\in I$, where $F(u)=|u|^\frac{8}{d-4}u$. We refer to
$I$ as the lifespan of $u$. We say that $u$ is a maximal-lifespan
solution if the solution cannot be extended to any strictly larger
interval. We say that $u$ is a global solution if $I=\Bbb R$.
\end{definition}
We define the scattering size of a solution to (\ref{fnls}) on a
time interval $I$ by
\begin{equation*}
S_I(u):=\int_I\int_{\Bbb R^d}|u(t,x)|^\frac{2(d+4)}{d-4}dxdt.
\end{equation*}
If $I=\Bbb R$, we write $S_{\Bbb R}(u)=S(u)$.

Associated to the notion of solution is a corresponding notion of
blowup, which precisely corresponds to the impossibility of
continuing the solution.
\begin{definition}[Blow up] We say that a solution $u$ to
(\ref{fnls}) blows up forward in time if there exists a time $t_1\in
I$ such that
\begin{equation*}
S_{[t_1, \sup(I))}(u)=\infty
\end{equation*}
and that $u$ blows up backward in time if there exists a time
$t_1\in I$ such that
\begin{equation*}
S_{(\inf{(I)}, t_1]}(u)=\infty.
\end{equation*}
\end{definition}

Sharp dispersive estimates for the biharmonic Schr\"odinger operator
in (\ref{gfnls}), namely for the linear group associated to
$i\partial_t +\Delta^2\pm\Delta$, have recently been obtained in
Ben-Artzi, Koch, and Saut \cite{BA}, while specific nonlinear fourth
order Schr\"odinger equations as in (\ref{gfnls}) have been recently
discussed in Fibich, Ilan, and Papanicolaou \cite{FIP}, Guo and Wang
\cite{GW}, Hao, Hsiao, and Wang \cite{HHW1}, \cite{HHW2}, Miao and
Zhang \cite{mz} and Segata \cite{Se3}. 
 In \cite{Pa}, B. Pausader established the global wellposedness in the defocusing
subcritical case, namely, $f(u)=|u|^{p-1}u$ with
$1<p<1+\frac{8}{d-4}$. Moreover, he estabilished the global
wellposedness and scattering for radial data in the defocusing
critical case, namely, $p=1+\frac{8}{d-4}$, where a very important
Strichartz estimate was established
\begin{lemma}[Strichartz estimates, \cite{Pa}] If $(q,r)$ is such that $\frac{2}{q}+\frac{d}{r}=\frac{d}{2}$,
where $2\leq q,r\leq\infty$ and $(q,r,d)\neq(2,\infty,2)$. Let $u$
be the solution of
\begin{equation} \label{inho}
\left\{ \aligned
    iu_t +  \Delta^2 u  & = h,\\
     u(t_0) \in & \ \dot{H}^2(\mathbb{R}^d).
\endaligned
\right.
\end{equation}
Then we have \begin{equation*} \|\Delta u\|_{L_t^qL_x^r(I\times\Bbb
R^d)}\lesssim \|\Delta u_0\|_2+\|\nabla
h\|_{L_t^2L_x^\frac{2d}{d+2}(I\times\Bbb R^d)}.
\end{equation*}
\end{lemma}
 The key feature of such lemma is that the
spacetime norm of the second derivative of $u$ is estimated using
only one derivative of the forcing term. In fact, this is the
consequence of smoothing effect for all higher order nonlinear
Schr\"odinger equation, see Proposition 2 in \cite{mz}. This is in
sharp contrast with the classical second order nonlinear
Schr\"odinger equations, where the estimate like (\ref{inho}) does
not hold true as it would violate Galilean invariance. Moreover,
local well-posedness and stability were established
\begin{theorem}[Local well-posedness, \cite{Pa}]\label{lw} Let
$d\geq5$. Given $u_0\in\dot{H}_x^2(\Bbb R^d)$ and $t_0\in \Bbb R$,
there exists a unique maximal-lifespan solution $u: I\times\Bbb
R^d\rightarrow\Bbb C$ to (\ref{fnls}) with initial data
$u(t_0)=u_0$. This solution also has the following properties:
\renewcommand{\labelenumi}{$\bullet$}
\begin{enumerate}
\item (Local existence) $I$ is an open neighborhood of $t_0$.
\item (Energy conservation) The energy of $u$ is conserved, that is,
$E(u(t))=E(u_0)$ for all $t\in I$.
\item (Continuous dependence) If $u_0^{(n)}$ is a sequence
converging to $u_0$ in $\dot{H}^2_x(\Bbb R^d)$ and $u^{(n)}:
I_n\times\Bbb R^d\rightarrow\Bbb C$ are the associated
maximal-lifespan solutions, then $u^{(n)}$ converge locally
uniformly to $u$, that is, on every  compact interval $K\subset I$,
and $K\subset I_n$ for all sufficiently large $n$, $u_n$ converges
strongly to $u$ in $C_t^0\dot{H}_x^2(K\times\Bbb R^d)\cap
L_{t,x}^\frac{2(d+4)}{d-4}(K\times\Bbb R^d)$ as
$n\rightarrow\infty$.
\item (Blowup criterion) If $\sup(I)$ is finite, then $u$ blows up
forward in time; if $\inf(I)$ is finite, then $u$ blows up backward
in time.
\item  (Scattering) If $\sup(I)=\infty$ and $u$ does not blow up forward in time, then $u$ scatters forward in time, that is, there
exists a unique $u_+\in \dot{H}_x^2(\Bbb R^d)$ such that
\begin{equation}\label{equ31}
\lim_{t\rightarrow\infty}\|u(t)-e^{it\Delta^2}u_+\|_{\dot{H}_x^2(\Bbb
R^d)}=0.
\end{equation}
Conversely, given $u_+\in \dot{H}_x^2(\Bbb R^d)$ there is a unique
solution to (\ref{fnls}) in a neighborhood of infinity so that
(\ref{equ31}) holds..
\item (Small data global existence) If $\|\Delta u_0\|_2$ is
sufficiently small (depending on $d$), then $u$ is a global solution
which does not blow up either forward or backward in time. Indeed,
in this case $S_{\Bbb R}(u)\lesssim \|\Delta
u_0\|_2^\frac{2(d+4)}{d-4}$.
\end{enumerate}
\end{theorem}
\begin{theorem}[Stability, \cite{Pa}]\label{stability}
Let $d\geq5$. Let $I\subset \Bbb R$ be a compact time interval such
that $0\in I$, and $\tilde{u}$ be an approximate solution of
(\ref{fnls}) in the sense that
\begin{equation*}
i\partial_t\tilde{u}+\Delta^2\tilde{u}-|\tilde{u}|^\frac{8}{d-4}\tilde{u}=e
\end{equation*}
for some $e$ with $\nabla e\in L_t^2L_x^\frac{2d}{d+2}(I\times\Bbb
R^d)$. Assume that
\begin{equation*}\|\tilde{u}\|_{L_{t,x}^\frac{2(d+4)}{d-4}(I\times\Bbb
R^d)}<+\infty\quad \text{and}\quad
\|\tilde{u}\|_{L_t^\infty\dot{H}^2_x(I\times\Bbb
R^d)}<+\infty.\end{equation*} For any $\Lambda>0$ there exists
$\delta_0>0$ such that if
\begin{equation*}
\|\nabla e\|_{L_t^2L_x^\frac{2d}{d+2}(I\times\Bbb R^d)}\leq\delta
\end{equation*}
and if $u_0\in \dot{H}^2_x(\Bbb R^d)$ satisfies
\begin{equation*}
\|\tilde{u}(0)-u_0\|_{\dot{H}^2}\leq \Lambda\ \  \text{and}\ \
\|\nabla
e^{it\Delta^2}(\tilde{u}(0)-u_0)\|_{L_t^\frac{2(d+4)}{d-4}L_x^\frac{2d(d+4)}{d^2-2d+8}(I\times\Bbb
R^d)}\leq\delta
\end{equation*}
for some $\delta\in(0, \delta_0)$, then there exists a solution
$u\in C(I, \dot{H}^2)$ of (\ref{fnls}) such that $u(0)=u_0$.
Moreover, $$\|u\|_{L_{t,x}^\frac{2(d+4)}{d-4}(I\times\Bbb
R^d)}<\infty.$$
\end{theorem}
Let
\begin{equation*}W(x)=W(x,t)=\Big(\frac{\big(d(d-4)(d^2-4)\big)^\frac{1}{4}}{1+|x|^2}\Big)^\frac{d-4}{2}
\end{equation*}
be a stationary solution of (\ref{fnls}). That is $W\geq0$ solves
the nonlinear elliptic equation
\begin{equation}\label{equ41}
\Delta^2 W=|W|^\frac{8}{d-4}W.
\end{equation}
Analogous to the nonlinear Schr\"odinger equations of the second
order, we have
\begin{conjecture}\label{cj} Let $d\geq5$ and let $u: I\times\Bbb R^d\rightarrow\Bbb
C$ be a solution to (\ref{fnls}) and $W$ is the stationary solution
of this equation. If
\begin{equation}\label{assumption} E_*:=\sup_{t\in I}\|\Delta
u(t)\|_{2}<\|\Delta W\|_{2},
\end{equation}
then \begin{equation*} \int_I\int_{\Bbb
R^d}|u(t,x)|^\frac{2(d+4)}{d-4}dxdt\leq C(E_*)<\infty.
\end{equation*}
\end{conjecture}

Naturally, we will apply the ideas and techniques which come from
the study of classical focusing nonlinear Schr\"odinger equations to
fourth order nonlinear Schr\"odinger equations. For energy critical
nonlinear Schr\"odinger equations
\begin{equation}\label{nls} \left\{ \aligned
    iu_t +  \Delta u  & = \lambda|u|^\frac{4}{d-2}u, \quad  \text{in}\  \mathbb{R}^d \times \mathbb{R},\\
     u(0)&=u_0(x)\in\dot{H}_x^1(\Bbb R^d),
\endaligned
\right.
\end{equation}
the local well-posedness and global well-posedness for small data
were established by T. Cazenave and F. B. Weissler \cite{CW}
regardless of the sign of $\lambda$. There have been a lot of works
devoted to obtaining the global well-posedness and scattering for
large data in defocusing case $\lambda=1$, see \cite{borg:scatter},
\cite{ckstt:gwp}, \cite{rv}, \cite{tao}, \cite{vis}, etc.

However, the global well-posedness and scattering for large data in
focusing case $\lambda=-1$ remains not completely solved . In
\cite{kenig-merle}, C. E. Kenig and F. Merle introduced an efficient
approach to deal with the focusing energy-critical nonlinear
Schr\"odinger equations, where they obtained global well-posedness
and scattering for radial data with energy and kinetic energy less
than those of ground state in the focusing case in dimensions $3\leq
d\leq5$. They reduced matters to a rigidity theorem using a
concentration compactness theorem, with the aid of localized Virial
identity. The radiality enters only at one point in the proof of the
rigidity theorem because of the difficulty in controlling the motion
of spatial center of global solutions. Moreover, their result is
sharp because the ground state itself is the solutions of
(\ref{nls}) but it does not scatter. One of the main ingredients in
their arguments is proved by S. Keraani in \cite{sh1}, namely, the
fact that every sequence of solutions to the linear Schr\"odinger
equation, with bounded data in $\dot{H}^1(\Bbb R^d)$ ($d\geq3$) can
be written, up to subsequence, as an almost orthogonal sum of
sequences of the type
$\lambda_n^{-\frac{d-2}{2}}V((t-t_n)/\lambda_n^2,
(x-x_n)/\lambda_n)$, where $V$ is a solution of the linear
Schr\"odinger equation, with a small remainder term in Strichartz
norm. Earlier steps in this direction include \cite{BG}.

R. Killip and M. Visan \cite{KM} extended C. E. Kenig and F. Merle's
result to nonradial case in $d\geq5$. The method is to reduce
minimal kinetic energy blow up solutions to almost periodic
solutions modulo symmetries, which match one of the three scenarios:
finite time blowup, low-to-high frequency cascade and soliton. Then
the aim is to eliminate such solutions. The finite time blowup
solutions can be precluded using the method in \cite{kenig-merle}.
For the other two types of solutions, R. Killip and M. Visan proved
that they admit additional regularities, namely, they belong to
$L_t^\infty\dot{H}^{-\epsilon}_x$ for some $\epsilon>0$. In
particular, they are in $L_x^2$. Similar ideas have appeared in
\cite{KTV2} and \cite{KVZ} in order to deal with mass-critical
nonlinear Schr\"odinger equations. But different from before, a
remarkable difficulty comes from the minimal kinetic energy blowup
solution because the kinetic energy, unlike the energy, is not
conserved. Related arguments (for the cubic NLS in three spatial
dimensions) appeared in \cite{kem2}. Now the low-to-high frequency
cascade can be precluded by negative regularity and the conservation
of mass. To preclude the soliton, one need to control of motion of
spatial center of the soliton solution. The method comes from
\cite{dhr} and \cite{kem} with the aid of negative regularity. The
fist step is to note that a minimal kinetic energy blowup solution
with finite mass must have zero momentum. A second ingredient is a
compactness property of the orbit of $\{u(t)\}$ in $L_x^2$. Finally
the soliton-like solution is precluded by using a truncated Virial
identity. Note that the negative regularity in \cite{KM} cannot be
obtained in dimensions three and four because the dispersion is too
weak. Indeed the method of \cite{kenig-merle} and \cite{KM} can be
applied to defocusing case without much difficulty.

In this paper, we will verify Conjecture \ref{cj} in radial case,
namely,
\begin{theorem}[Spacetime bounds]\label{tm1} Let $d\geq5$ and let
 $u: I\times\Bbb R^d\rightarrow\Bbb
C$ be a radial solution to (\ref{fnls}). If
\begin{equation*}\label{equ01} E_*:=\sup_{t\in I}\|\Delta
u(t)\|_{2}<\|\Delta W\|_{2},
\end{equation*}
then \begin{equation*} \int_I\int_{\Bbb
R^d}|u(t,x)|^\frac{2(d+4)}{d-4}dxdt\leq C(E_*)<\infty.
\end{equation*}
\end{theorem}
A more effective criterion for global well-posedness (depending
directly on $u_0$) can be obtained using an energy-trapping argument
in Section 3 (the corresponding argument for nonlinear Schr\"odinger
equations is in \cite{kenig-merle}).
\begin{corollary} Let $d\geq5$ and let $u_0\in \dot{H}_x^2(\Bbb
R^d)$ be a radial function and such that $\|\Delta u_0\|_2<\|\Delta
W\|_2$ and $E(u_0)<E(W)$. Then the corresponding solution $u$ to
(\ref{fnls}) is global and moreover,
\begin{equation*}
\int_{\Bbb R}\int_{\Bbb R^d}|u(t,x)|^\frac{2(d+4)}{d-4}dtdx<\infty.
\end{equation*}
\end{corollary}

In this paper, we establish the corresponding result of the theorem
in \cite{kenig-merle} on the setting of nonlinear Schr\"odinger
equations of fourth order. For later use in \cite{mxz3}, the
arguments here are direct ``fourth order'' analogue of \cite{KM},
including \cite{kenig-merle}, \cite{sh1} and \cite{KTV2}. First, we
will do a lot of ground work including establishing concentration
compactness principle and the energy-trapping of the ground state.
Next, we reduce the failure of Conjecture \ref{cj} to almost
periodic solutions, where we will rely heavily on Theorem
\ref{stability}. To show that such almost periodic solutions match
one of the three scenarios, we analyze the properties of the almost
periodic solutions such as quasi-uniqueness of $N$, compactness of
almost periodic solutions, etc (see Section 4). Because we are
considering the minimal kinetic blow up solution, the assumption
(\ref{assumption}) plays an important role, which is used in the
proof of quasi-uniqueness of $N$. Finally, we established localized
Virial identity and precluded all the possibility of the three
scenarios. Note that the radiality enters only in Section 8, so all
the conclusions in Section 3- Section 7 remain true for general
solutions. Moreover, the method here applies equally well to
defocusing nonlinear Schr\"odinger equations of fourth order.
Because no Galilean transformation is available for (\ref{fnls}), it
seems difficult to remove the radial assumption as in \cite{KM} even
in high dimensions. But we can remove the radial assumption in the
defocusing case in dimensions $d\geq9$, see \cite{mxz3}.

After the paper was finished and submitted, we learned that B.
Pausader has obtained independently the similar result in
\cite{pa2}, where the author proved the same result with
(\ref{assumption}) replaced by $E(u)<E(W)$.

The rest of the paper is organized as follows: In Section 2, we
introduce some notations. The energy-trapping of the ground state is
given in Section 3. In Section 4, we define almost periodic
solutions and list their properties. The concentration compactness
principle is proved in Section 5. In Section 6, we reduce the
failure of Conjecture \ref{cj} to the existence of almost periodic
solutions and in Section 7, we prove that such solutions must admit
one of three scenarios, namely, we set up three enemmies. Finally,
we preclude all the scenarios in Section 8.
\section{Notations}
 We introduce some notations. If $X, Y$ are
nonnegative quantities, we use $X\lesssim Y $ or $X=O(Y)$ to denote
the estimate $X\leq CY$ for some $C$ which may depend on the energy
$E(u)$ and $X \sim Y$ to denote the estimate $X\lesssim Y\lesssim
X$. Sometimes we write $X\sim_{c,C,u}Y$ to mean the implicit
constant depends on $c$, $C$ and $E(u)$. We use $X\ll Y$ to mean $X
\leq c Y$ for some small constant $c$ which is again allowed to
depend on $E(u)$. We write $L_t^qL_x^r$ to denote the Banach space
with norm
$$\|u\|_{L_t^qL_x^r(\mathbb{R}\times\mathbb{R}^d)}:=\Big(
\int_{\mathbb{R}}\big(\int_{\mathbb{R}^d}|u(t,x)|^rdx\big)^{q/r}dt\Big)^{1/q},$$
with the usual modifications when $q$ or $r$ are equal to infinity,
or when the domain $\mathbb{R}\times\mathbb{R}^d$ is replaced by
spacetime slab such as $I\times\mathbb{R}^d$. When $q=r$ we
abbreviate $L_t^qL_x^q$ as $L_{t,x}^q$.

We use $C\gg1$ to denote various large finite constants. and $0< c
\ll 1$ to denote various small constants.

The Fourier transform on $\mathbb{R}^d$ is defined by
\begin{equation*}
\aligned \widehat{f}(\xi):= \big( 2\pi
\big)^{-\frac{d}{2}}\int_{\mathbb{R}^d}e^{- ix\cdot \xi}f(x)dx ,
\endaligned
\end{equation*}
giving rise to the fractional differentiation operators
$|\nabla|^{s}$,  defined by
\begin{equation*}
\aligned \widehat{|\nabla|^sf}(\xi):=|\xi|^s\widehat{f}(\xi).
\endaligned
\end{equation*}
These define the homogeneous Sobolev norms
\begin{equation*}
\big\|f\big\|_{\dot{H}^s_x}:= \big\| |\nabla|^sf
\big\|_{L^2_x(\mathbb{R}^d)}.
\end{equation*}

Let $e^{it\Delta^2 }$ be the free fourth order Schr\"{o}dinger
propagator given by
\begin{equation*}
\widehat{e^{it\Delta^2 } f}(\xi)=e^{it|\xi|^4}\widehat{f}(\xi).
\end{equation*}
We recall some basic facts in Littlewood-Paley theory. Let
$\varphi(\xi)$ be a radial bump function supported in the ball
$\{\xi\in\mathbb{R}^d: |\xi|\leq\frac{11}{10}\}$ and equal to 1 on
the ball $\{\xi\in\mathbb{R}^d: |\xi|\leq1\}$. For each number
$N>0$, we define the Fourier multipliers
\begin{align*}
\widehat{P_{\leq N}f}(\xi)&:=\varphi(\xi/N)\hat{f}(\xi),\\
\widehat{P_{\geq N}f}(\xi)&:=(1-\varphi(\xi/N))\hat{f}(\xi),\\
\widehat{P_{N}f}(\xi)&:=(\varphi(\xi/N)-\varphi(2\xi/N))\hat{f}(\xi)
\end{align*}
and similarly $P_{<N}$ and $P_{\geq N}$. We also define
$$P_{M<\cdot\leq N}:=P_{\leq N}-P_{\leq M}=\sum_{M<N'\leq N}P_{N'}$$
whenever $M<N$. We will usually use these multipliers when $M$ and
$N$ are dyadic numbers; in particular, all summations over $N$ or
$M$ are understood to be over dyadic numbers. Nevertheless, it will
occasionally be convenient to allow $M$ and $N$ to not be a power of
2. Note that $P_N$ is not truly a projection; to get around this, we
will occasionally need to use fattened Littlewood-Paley operators:
\begin{equation}
\tilde{P}_N:=P_{N/2}+P_N+P_{2N}.
\end{equation}
They obey $P_N\tilde{P}_N=\tilde{P}_NP_N=P_N$.

As all Fourier multipliers, the Littlewood-Paley operators commute
with the propagator $e^{it\Delta^2}$, as well as with differential
operators such as $i\partial_t+\Delta^2$. We will use basic
properties of these operators many times, including
\begin{lemma}[Bernstein estimates]\label{Bern} For $1\leq p\leq q\leq\infty$,
\begin{align*}\||\nabla|^{\pm s}P_Nf\|_{L_x^p(\mathbb{R}^d)}&\sim
N^{\pm s}\|P_Nf\|_{L_x^p(\mathbb{R}^d)},\\
\|P_{\leq N}f\|_{L_x^q(\mathbb{R}^d)}&\lesssim
N^{\frac{d}{p}-\frac{d}{q}}\|P_{\leq N}f\|_{L_x^p(\mathbb{R}^d)},\\
\|P_{N}f\|_{L_x^q(\mathbb{R}^d)}&\lesssim
N^{\frac{d}{p}-\frac{d}{q}}\|P_{N}f\|_{L_x^p(\mathbb{R}^d)}.
\end{align*}
\end{lemma}
\section{Ground state}
Let
\begin{equation*}W(x)=W(x,t)=\Big(\frac{\big(d(d-4)(d^2-4)\big)^\frac{1}{4}}{1+|x|^2}\Big)^\frac{d-4}{2}
\end{equation*}
be a stationary solution of (\ref{fnls}). That is $W\geq0$ solves
the nonlinear elliptic equation
\begin{equation}\label{equ41}
\Delta^2 W=|W|^\frac{8}{d-4}W,
\end{equation}
then by the invariances of the equation,  for
$\theta_0\in[-\pi,\pi]$, $\lambda_0>0$, $x_0\in \Bbb R^d$,
\begin{equation*}
W_{\theta_0,x_0,\lambda_0}=e^{i\theta_0}{\lambda_0}^\frac{d-4}{2}W(\lambda_0(x-x_0))
\end{equation*}
is still a solution. By the work of \cite{Ed}, we have the following
characterization of $W$:
\begin{equation}\label{equ46}
\forall u\in \dot{H}^2, \ \ \|u\|_{L^{2^{\#}}}\leq C_d\|\Delta
u\|_{L^2};
\end{equation}
moreover, if $u\neq0$ is such that
\begin{equation}\label{equ42}\|u\|_{L^{2^\#}}= C_d\|\Delta u\|_{L^2},\end{equation} then
there exist $(\theta_0, \lambda_0, x_0)$  such that
$u=W_{\theta_0,x_0,\lambda_0}$, where $C_d$ is the best constant of
the Sobolev inequality in dimension $d$ and $2^{\#}=\frac{2d}{d-4}$.

The equation (\ref{equ41}) gives $\int|\Delta
W|^2=\int|W|^{2^{\#}}$. Also, (\ref{equ42}) yields $C_d^2\int|\Delta
W|^2=\Big(\int|W|^{2^{\#}}\Big)^\frac{d-4}{d}$, an easy computation
shows that
\begin{equation*}
\int|\Delta W|^2=C_d^{-d/2}\ \ \text{and}\ \
E(W)=\frac{2}{d}C_d^{-d/2}.
\end{equation*}
\begin{lemma}\label{lem41}Assume that
\begin{equation*}
\|\Delta u\|_2<\|\Delta W\|_2
\end{equation*}
and $E(u)\leq (1-\delta_0)E(W)$ where $\delta_0>0$. Then there
exists $\bar{\delta}=\bar{\delta}(\delta_0, d)$ such that
\begin{equation}\label{equ43}
\int|\Delta u|^2\leq (1-\bar{\delta})\int|\Delta W|^2;
\end{equation}
\begin{equation}\label{equ44}
\int(|\Delta u|^2-|u|^{2^{\#}})\geq\bar{\delta}\int|\Delta u|^2;
\end{equation}
\begin{equation}\label{equ45}
E(u)\geq0.
\end{equation}
\end{lemma}
{\it Proof.} Consider the function
$f(y)=\frac{1}{2}y-\frac{d-4}{2d}C_d^{2^{\#}}y^\frac{2^\#}{2}$ and
let $\bar{y}=\|\Delta u\|_2^2$. From (\ref{equ46}), we have
\begin{equation*}
f(\bar{y})\leq E(u)\leq
(1-\delta_0)E(W)=(1-\delta_0)\frac{2}{dC_d^{d/2}}.
\end{equation*}
Note that $f(0)=0$,
$f'(y)=\frac{1}{2}-\frac{C_d^{2^{\#}}}{2}y^{\frac{2^{\#}}{2}-1}$, so
$f'(y)=0$ if and only if  $y=y_C$, where $y_C=C_d^{-d/2}=\int|\Delta
W|^2$. Note also that $f(y_C)=E(W)$. But since $0<\bar{y}<y_C$ and
$f(\bar{y})\leq (1-\delta_0)f(y_C)$ and $f$ is nonnegative and
strictly increasing between $0$ and $y_C$, $f''(y_C)\neq0$, we have
$0<f(\bar{y})$ and $\bar{y}\leq(1-\bar{\delta})\int|\Delta W|^2$.
Thus (\ref{equ43}) and (\ref{equ45}) hold.

To show (\ref{equ44}), consider the function
$g(y)=y-C_d^{2^{\#}}y^\frac{d}{d-4}$. Because of (\ref{equ46}), we
have that $\int(|\Delta u|^2-|u|^{2^{\#}})\geq g(\bar{y})$. Note
that $g(y)=0$ if and only if $y=0$ or $y=y_C$ and that $g'(0)=1$,
$g'(y_C)=-\frac{4}{d-4}$. We then have, for $0<y<y_C$, $g(y)\geq
C\min\{y, y_C-y\}$, so (\ref{equ44}) follows from $0\leq
\bar{y}<(1-\bar{\delta})y_C$ which is given by (\ref{equ43}).

By energy conservation, Lemma \ref{lem41} and a continuity argument,
we have
\begin{theorem}[Energy trapping]\label{lem42}
Let $u$ be a solution of (\ref{fnls}) with initial data $u_0$ such
that
\begin{equation*}
\int|\Delta u_0|^2<\int|\Delta W|^2\ \ \text{and}\ \
E(u_0)<(1-\delta_0)E(W).
\end{equation*}
Let $I\ni0$ be the maximal interval of existence. Let
$\bar{\delta}=\bar{\delta}(\delta_0, d)$ be as in Lemma \ref{lem41},
then for each $t\in I$, we have
\begin{equation}\label{equ403}
\int|\Delta u(t)|^2\leq (1-\bar{\delta})\int|\Delta W|^2;
\end{equation}
\begin{equation}\label{equ404}
\int(|\Delta u(t)|^2-|u(t)|^{2^{\#}})\geq\bar{\delta}\int|\Delta
u(t)|^2;
\end{equation}
\begin{equation}\label{equ405}
E(u(t))\geq0.
\end{equation}
\end{theorem}
{\it Proof.} See \cite{kenig-merle}.
\begin{corollary}\label{co41}
Let $u$, $u_0$ be as in Theorem \ref{lem42}. Then for all $t\in I$
we have $E(u(t))\simeq\int|\Delta u(t)|^2\simeq \int|\Delta u_0|^2$,
with comparability constants which depend only on $\delta_0$.
\end{corollary}
{\it Proof.} $E(u(t))\leq\int|\Delta u|^2dx$. But by (\ref{equ404})
we have
\begin{align*}
E(u)\geq&\ (\frac{1}{2}-\frac{1}{2^{\#}})\int|\Delta
u(t)|^2dx+\frac{1}{2^{\#}}\int\big(|\Delta
u(t)|^2-|u(t)|^{2^{\#}}\big)dx\\
\geq&\ C_{\bar{\delta}}\int|\Delta u(t)|^2dx,
\end{align*}
so the first equivalence follows. For the second one, note that
$E(u(t))=E(u_0)\simeq\int|\Delta u_0|^2dx$, by the first equivalence
when $t=0$.

\section{Almost periodic solutions}
\begin{definition}[Symmetry group]\label{d2}
For any phase $\theta\in \Bbb R/2\pi\Bbb Z$, position $x_0\in\Bbb
R^d$ and scaling parameter $\lambda>0$, we define the unitary
transformation $g_{\theta, x_0, \lambda}: \dot{H}^2(\Bbb
R^d)\rightarrow\dot{H}^2(\Bbb R^d)$ by the formula
\begin{equation*}[g_{\theta, x_0,
\lambda}f](x):=\lambda^{-\frac{d-4}{2}}e^{i\theta}f\big(\lambda^{-1}(x-x_0)\big).
\end{equation*}
We let $G$ be the collection of such transformations. If $u:
I\times\Bbb R^d\rightarrow\Bbb C$ is a function, we define
$T_{g_{\theta, x_0, \lambda}}u: \lambda^4I\times\Bbb
R^d\rightarrow\Bbb C$ where $\lambda^4I:=\{\lambda^4t: t\in I\}$ by
the formula
\begin{equation*}
[T_{g_{\theta, x_0,
\lambda}}u](t,x):=\lambda^{-\frac{d-4}{2}}e^{i\theta}u\big(\lambda^{-4}t,
\lambda^{-1}(x-x_0)\big).
\end{equation*}
We also let $G_{rad}\subset G$ denote the collection of
transformations in $G$ which preserve spherical symmetry, or more
explicitly,
\begin{equation*}
G_{rad}:=\{g_{\theta, 0, \lambda}: \theta\in \Bbb R/2\pi\Bbb Z;
\lambda>0\}.
\end{equation*}
\end{definition}

\begin{definition}[Enlarged group] For any $\theta\in \Bbb R/2\pi\Bbb Z$, position $x_0\in \Bbb R^d$, scaling parameter $\lambda>0$
 and time $t_0$, we define the unitary transformation $g_{\theta, x_0, \lambda, t_0}: \dot{H}^2_x(\Bbb R^d)\rightarrow\dot{H}^2_x(\Bbb R^d)$
 by the formula
 \begin{equation*}
g_{\theta, x_0, \lambda, t_0}=g_{\theta, x_0,
\lambda}e^{it_0\Delta^2}.
 \end{equation*}
Let $G'$ be the collection of such transformations. We also let $G'$
act on global spacetime functions $u: \Bbb R\times\Bbb
R^d\rightarrow\Bbb C$ by defining
\begin{equation*}
\big(T_{g_{\theta, x_0, \lambda,
t_0}}u\big)(t,x):=\frac{1}{\lambda^\frac{d-4}{2}}e^{i\theta}\big(e^{it_0\Delta^2}u\big)\big(\frac{t}{\lambda^4},
\frac{x-x_0}{\lambda}\big)
\end{equation*} Given any two sequences $g_n$, $g_n'$ in $G'$, we say that $g_n$
and $g_n'$ are asymptotically orthogonal if $(g_n)^{-1}g_n'$
diverges to infinity in $G$. If we write explicitly
\begin{equation*}
g_n=g_{\theta_n, x_n, \lambda_n, t_n}, \ \ g_n'=g_{\theta_n', x_n',
\lambda_n', t_n'},
\end{equation*}
then the asymptotic orthogonality is equivalent to
\begin{equation*}
\lim_{n\rightarrow\infty}\Big(\frac{\lambda_n}{\lambda_n'}+\frac{\lambda_n'}{\lambda_n}
+\frac{|t_n\lambda_n^4-t_n'(\lambda_n')^4|}{\lambda_n^2\lambda_n'^2}+\frac{|x_n-x_n'|^2}{\lambda_n\lambda_n'}\Big)=+\infty.
\end{equation*}
\end{definition}
\begin{remark}\label{rem1}If $g_n$, $g_n'\in G'$ are asymptotically orthogonal, then
\begin{equation}\label{equ512}
\lim_{n\rightarrow\infty}\langle \Delta g_nf, \Delta
g_n'f'\rangle_{L_x^2(\Bbb R^d)}=0,\ \ \text{for all}\ \ f, f'\in
\dot{H}_x^2(\Bbb R^d).
\end{equation}
A variant of this is that if $v$, $v'\in
L_{t,x}^\frac{2(d+4)}{d-4}(\Bbb R\times\Bbb R^d)$, then
\begin{equation}\label{equ511}
\lim_{n\rightarrow\infty}\big\||T_{g_n}v|^{1-\theta}|T_{g_n'}v'|^\theta\big\|_{L_t^qL_x^r}=0
\end{equation}
for any $0<\theta<1$ and admissible pair $(q, r)$ ($q<\infty$), that
is, $\frac{4}{q}+\frac{d}{r}=\frac{d}{2}-2$.
\end{remark}
\begin{definition}[Almost periodic solutions]\label{d1} Let $d\geq 5$. A solution $u$ to (\ref{fnls})
with lifespan $I$ is said to be almost periodic modulo $G$ if there
exist functions $N: I\rightarrow\Bbb R^+$, $x: I\rightarrow\Bbb R^d$
and $C: \Bbb R^+\rightarrow\Bbb R^+$ such that for all $t\in I$, and
$\eta>0$,
\begin{equation}\label{equ51}
\int_{|x-x(t)|\geq C(\eta)/N(t)}|\Delta u(t,x)|^2dx\leq\eta
\end{equation}
and
\begin{equation}\label{equ52}
\int_{|\xi|\geq C(\eta)N(t)}|\xi|^4|\hat{u}(t,\xi)|^2d\xi\leq\eta.
\end{equation}
We refer to the function $N$ as the frequency scale function for the
solution $u$, $x$ the spatial center function, and to $C$ as the
compactness modulus function.
\end{definition}
By Ascoli-Arzela theorem, almost periodicity modulo $G$ means that
the quotient orbit $\Big\{Gu(t): t\in I\Big\}$ is a precompact set
of $G\backslash \dot{H}^2$, where $G\backslash \dot{H}^2$ is the
moduli space of $G$-orbits $Gf:=\{gf: g\in G\}$ of $\dot{H}^2(\Bbb
R^d)$. Moreover, a family of functions is precompact in
$\dot{H}_x^2$ if and only if it is norm-bounded and there exists a
compactness modulus function $C$ so that
\begin{equation*}
\int_{|x|\geq C(\eta)}|\Delta f(x)|^2dx+\int_{|\xi|\geq
C(\eta)}|\xi|^4|\hat{f}(\xi)|^2d\xi\leq \eta
\end{equation*}
for all functions $f$ in the family. By Sobolev embedding, any
solution $u: I\times\Bbb R^d\rightarrow\Bbb C$ that is almost
periodic modulo $G$ must also satisfy \begin{equation*}
\int_{|x-x(t)|\geq C(\eta)/N(t)}|u(t,x)|^\frac{2d}{d-4}dx\leq\eta.
\end{equation*}
\begin{lemma}[Quasi-uniqueness of $N$]\label{lem55}Let $u$ be a
non-zero solution to (\ref{fnls}) satisfying (\ref{assumption}) with
lifespan $I$ that is almost periodic modulo $G$ with frequency scale
function $N: I\rightarrow\Bbb R^+$ and compactness modulus function
$C: \Bbb R^+\rightarrow\Bbb R^+$ and also almost periodic modulo $G$
with frequency scale function $N': I\rightarrow\Bbb R^+$ and
compactness modulus function $C': \Bbb R^+\rightarrow\Bbb R^+$. Then
we have
\begin{equation*}
N(t)\sim_{u, C, C'}N'(t)
\end{equation*}
for all $t\in I$.
\end{lemma}
{\it Proof.} It suffices to prove $N'(t)\lesssim_{u, C, C'}N(t)$,
for all $t\in I$. Otherwise, there exists a sequence $\{t_n\}$ such
that $\lim_{n\rightarrow\infty}N(t_n)/N'(t_n)=0$. For any $\eta>0$,
by Definition \ref{d1}, we have
\begin{equation*}
\int_{|x-x'(t_n)|\geq C'(\eta)/N'(t_n)}|\Delta u(t_n,x)|^2dx\leq\eta
\end{equation*}
and
\begin{equation}\label{equ411}
\int_{|\xi|\geq
C(\eta)N(t_n)}|\xi|^4||\hat{u}(t_n,\xi)|^2d\xi\leq\eta.
\end{equation}
Let $u(t_n,x)=u_1(t_n,x)+u_2(t_n,x)$, where
$u_1(t_n,x)=u(t_n,x)1_{|x-x'(t_n)|\geq C'(\eta)/N'(t_n)}$,
$u_2(t_n,x)=u(t_n,x)1_{|x-x'(t_n)|< C'(\eta)/N'(t_n)}$. Then by
Plancherel's theorem, we have
\begin{equation}\label{equ53}
\int_{\Bbb R^d}|\xi|^4|\hat{u}_1(t_n,\xi)|^2d\xi\lesssim \eta,
\end{equation}
while from Cauchy-Schwartz we have \begin{equation*}
\sup_{\xi\in\Bbb R^d}|\xi|^4|\hat{u}_2(t_n,\xi)|^2\lesssim_{\eta,
C'}\|\Delta u(t_n)\|_{2}^2N'(t_n)^{-d}. \end{equation*} Integrating
the last inequality over the ball $|\xi|\leq C(\eta)N(t_n)$, we get
\begin{align*}
\int_{|\xi|\leq C(\eta)N(t)}|\xi|^4|\hat{u}(t_n,\xi)|^2d\xi&\
\lesssim\int_{\Bbb
R^d}|\xi|^4|\hat{u}_1(t_n,\xi)|^2d\xi+\int_{|\xi|\leq
C(\eta)N(t)}|\xi|^4|\hat{u}_2(t_n,\xi)|^2d\xi\\
&\ \lesssim\eta+O_{\eta, C, C'}(\|\Delta
u(t_n)\|_2^2N(t_n)^dN'(t_n)^{-d}).
\end{align*}
This, combined with (\ref{equ411}), (\ref{assumption}) and Corollary
\ref{co41}, yields that
\begin{align*}
\int|\Delta u_0(x)|^2dx\sim\int_{\Bbb
R^d}|\xi|^4|\hat{u}(t_n,\xi)|^2d\xi \lesssim&\ \eta+O_{\eta, C,
C'}\big(\|\Delta u(t_n)\|_2^2N(t_n)^dN'(t_n)^{-d}\big)\\
\lesssim&\ \eta+O_{\eta, C, C'}\big(\|\Delta
W\|_2^2N(t_n)^dN'(t_n)^{-d}\big).\end{align*} Since
$\lim_{n\rightarrow\infty}N(t_n)/N'(t_n)=0$, we have
\begin{equation*}\int|\Delta u_0(x)|^2dx\lesssim\eta.
\end{equation*}
By the arbitrary of $\eta$, we get \begin{equation*}\int|\Delta
u_0(x)|^2dx=0.
\end{equation*}
Thus, $u_0\equiv0$ and by mass conservation, $u(t)\equiv0$  for all
$t\in\Bbb R$. This contradicts that $u$ is non-zero.
\begin{lemma}[Quasi-continuous dependence of $N$ on $u$]\label{lem56}
Let $u^{(n)}$ be a sequence of solutions to (\ref{fnls}) with
lifespan $I^{(n)}$ satisfying (\ref{assumption}), which are almost
periodic modulo scaling with frequency scale functions $N^{(n)}:
I^{(n)}\rightarrow\Bbb R^+$ and compactness modulus function $C:
\Bbb R^+\rightarrow\Bbb R^+$, independent of $n$. Suppose that
$u^{(n)}$ converge locally uniformly to a non-zero solution $u$ to
(\ref{fnls}) with lifespan $I$. Then $u$ is almost periodic modulo
scaling with a frequency scale function $N: I\rightarrow\Bbb R^+$
and compactness modulus function $C$. Furthermore, we have
\begin{equation}\label{equ54}
N(t)\sim_{u, C}\liminf_{n\rightarrow\infty}N^{(n)}(t)\sim_{u,
C}\limsup_{n\rightarrow\infty}N^{(n)}(t)
\end{equation}
for all $t\in I$. Finally, if all $u^{(n)}$ are spherically
symmetric, then $u$ is also.
\end{lemma}
{\it Proof.} We first show that
\begin{equation}\label{equ32}
0<\liminf_{n\rightarrow\infty}N^{(n)}(t)\leq\limsup_{n\rightarrow\infty}N^{(n)}(t)<+\infty
\end{equation}
for all $t\in I$. Indeed, if one of these inequalities failed for
some $t$, then (by passing to a subsequence if necessary)
$N^{(n)}(t)$ would converge to zero or infinity as
$n\rightarrow\infty$. Thus by Definition \ref{d1}, $u^{(n)}(t)$
would converge weakly to zero, hence by the local uniform
convergence, would converge strongly to zero. But this contradicts
the hypothesis that $u$ is not identically zero. This establishes
(\ref{equ32}).

From (\ref{equ32}), we see that for each $t\in I$ the sequence
$N^{(n)}(t)$ has at least one limit point $N(t)$. Thus, using the
local uniform convergence we easily verify that $u$ is almost
periodic modulo scaling with frequency scale function $N$ and
compactness modulus function $C$. It is also clear that if all
$u^{(n)}$ are spherically symmetric, then $u$ is also.

It remains to establish (\ref{equ54}), which we prove by
contradiction. Suppose it fails. Then given any $A=A_u$, there
exists a $t\in I$ for which $N^{(n)}(t)$ has at least two limit
points which are separated by a ratio of at least $A$, and so $u$
has two frequency scale functions with compactness modulus function
$C$ which are separated by this ratio. But this contradicts Lemma
\ref{lem55} for $A$ large enough depending on $u$. Hence
(\ref{equ54}) holds.
\begin{definition}[Normalised solution] Let $u$ be a solution to
(\ref{fnls}), which is almost periodic modulo $G$ with frequency
scale function $N$, position center function $x$. We say that $u$ is
normalised if the lifespan $I$ contains zero and
\begin{equation*}
N(0)=1,\ \ x(0)=0.
\end{equation*}
More generally, we can define the normalisation of a solution $u$ at
time $t_0$ in its lifespan $I$ to be
\begin{equation}\label{norm}
u^{[t_0]}:=T_{g_{0, -x(t_0)N(t_0), N(t_0)}}\big(u(\cdot+t_0)\big).
\end{equation}
Observe that $u^{[t_0]}$ is  a normalised solution which is almost
periodic modulo $G$ and has lifespan
\begin{equation*}
I^{[t_0]}:=\{s\in\Bbb R: t_0+s/N(t_0)^4\in I\}.
\end{equation*}
It has frequency scale function \begin{equation*}
N_{u^{[t_0]}}(t)=\frac{N(t_0+tN(t_0)^{-4})}{N(t_0)}
\end{equation*}
and position center function
\begin{equation*}
x_{u^{[t_0]}}(t)=N(t_0)[x(t_0+tN(t_0)^{-4})-x(t_0)].
\end{equation*}
\end{definition}
\begin{lemma}[Compactness of almost periodic solutions] \label{lem57} Let
$u^{(n)}$ be a sequence of normalised maximal-lifespan solutions to
(\ref{fnls}) satisfying (\ref{assumption}) with lifespan
$I^{(n)}\ni0$, which are almost periodic modulo $G$ with frequency
scale functions $N^{(n)}: I^{(n)}\rightarrow\Bbb R^+$ and a uniform
compactness modulus function $C: \Bbb R^+\rightarrow\Bbb R^+$.
Assume that we also have a uniform energy bound
\begin{equation*}0<\inf_nE(u^{(n)})\leq\sup_nE(u^{(n)})<\infty.\end{equation*}
Then after passing to a subsequence if necessary, there exists a
non-zero maximal-lifespan solution $u$ to (\ref{fnls}) with lifespan
$I\ni0$ that is almost periodic modulo $G$, such that $u^{(n)}$
converge locally uniformly to $u$. Moreover, if  all $u^{(n)}$ are
spherically symmetric and almost periodic modulo $G_{rad}$, then $u$
is also.
\end{lemma}
{\it Proof.} By hypothesis and Definition \ref{d1}, we see that for
every $\varepsilon>0$ there exists $R>0$ such that
\begin{equation*}
\int_{|x|\geq R}|\Delta u^{(n)}(0,x)|^2dx\leq\varepsilon
\end{equation*}
and
\begin{equation*}
\int_{|\xi|\geq
R}|\xi|^4|\widehat{u^{(n)}}(0,\xi)|^2d\xi\leq\varepsilon
\end{equation*}
for all $n$. Since $\sup_nE(u^{(n)})<\infty$, we have
$\sup_n\|\Delta u^{(n)}(0)\|_2^2<\infty$. By the Ascoli-Arzela
Theorem, we see that the sequence $u^{(n)}(0)$ is precompact in the
strong topology of $\dot{H}^2(\Bbb R^d)$. Thus passing to a
subsequence if necessary, we can find $u_0\in \dot{H}^2(\Bbb R^d)$
such that $u^{(n)}(0)$ converge strongly to $u_0$ in $\dot{H}^2(\Bbb
R^d)$. Since $0<\inf_nE(u^{(n)})$, we see that $u_0$ is not
identically zero.

Now let $u$ be the maximal Cauchy development of $u_0$ from time 0,
with lifespan $I$. By Theorem \ref{lw}, $u^{(n)}$ converges locally
uniformly to $u$. The remaining claims now follow from Lemma
\ref{lem56}.
\begin{corollary}[Local constancy of $N$]\label{c58}
Let $u$ be a non-zero maximal-lifespan solution to (\ref{fnls})
satisfying (\ref{assumption}) with lifespan $I$ that is almost
periodic modulo $G$ with frequency scale function $N:
I\rightarrow\Bbb R^+$. Then there exists a small number $\delta$,
depending on $u$, such that for every $t_0\in I$ we have
\begin{equation}\label{equ55}
[t_0-\delta N(t_0)^{-4},t_0+\delta N(t_0)^{-4}]\subset I
\end{equation}
and \begin{equation}\label{equ56} N(t)\sim_u N(t_0)
\end{equation}
whenever $|t-t_0|\leq\delta N(t_0)^{-4}$.
\end{corollary}
{\it Proof.} Let us establish (\ref{equ55}) first. We argue by
contradiction. Assume that (\ref{equ55}) failed. Then there exist
sequences $t_n\in I$ and $\delta_n\rightarrow0$ such that
$t_n+\delta_nN(t_n)^{-4}\not\in I$ for all $n$. Define the
normalisation $u^{[t_n]}$ of $u$ from time $t_n$ by (\ref{equ54}).
Then $u^{[t_n]}$ are maximal-lifespan normalised solutions whose
lifespan $I^{[t_n]}$ contain $0$ but not $\delta_n$; they are also
almost periodic modulo $G$ with frequency scale functions
\begin{equation*}
N^{[t_n]}(s):=N(t_n+sN(t_n)^{-4})/N(t_n)
\end{equation*}
and the same compactness modulus function $C$ as $u$. Applying Lemma
\ref{lem57} (and passing to a subsequence if necessary), we conclude
that $u^{[t_n]}$ converge locally uniformly to a maximal-lifespan
solution $v$ with some lifespan $J\ni0$. By Theorem \ref{lw}, $J$ is
open and so contains $\delta_n$ for all sufficiently large $n$. This
contradicts the local uniform convergence as, by hypothesis,
$\delta_n$ does not belong to $I^{[t_n]}$. Hence (\ref{equ55})
holds.

We now show (\ref{equ56}). Again, we argue by contradiction,
shrinking $\delta$ if necessary. Assume (\ref{equ56}) failed no
matter how small one select $\delta$. Then one can find sequences
$t_n$, $t_n'\in I$ such that $s_n:=(t_n'-t_n)N(t_n)^4\rightarrow0$
but $N(t_n')/N(t_n)$ converge to either zero or infinity. If we
define $u^{[t_n]}$ and $N^{[t_n]}$ as before and apply Lemma
\ref{lem57}, we see once again that $u^{[t_n]}$ converge locally
uniformly to maximal solution with some open interval $J\ni0$. But
then $N^{[t_n]}(s_n)$ converge to either zero or infinity and thus
by Definition \ref{d1}, $u^{[t_n]}(s_n)$ are converging weakly to
zero. On the other hand, since $s_n$ converge to zero and
$u^{[t_n]}$ are locally uniformly convergent to $v\in C_{t,
loc}^0\dot{H}^2_x(J\times\Bbb R^d)$, we may conclude that
$u^{[t_n]}(s_n)$ converge strongly to $v(0)$ in $\dot{H}_x^2(\Bbb
R^d)$. Thus $\|v(0)\|_{\dot{H}_x^2(\Bbb R^d)}=0$. So $v(0)\equiv0$.
Since $E(u^{[t_n]})=E(u)$, we see that $u$ vanishes. Thus
(\ref{equ56}) holds.

As a direct consequence of Corollary \ref{c58}, we have
\begin{corollary}[Blowup criterion]\label{c59} Let $u$ be a non-zero
maximal-lifespan solution to (\ref{fnls}) satisfying
(\ref{assumption}) with lifespan $I$ that is almost periodic modulo
$G$ with frequency scale function $N: I\rightarrow\Bbb R^+$. If $T$
is any finite endpoint of $I$, then
\begin{equation*}
\lim_{t\rightarrow T}N(t)=\infty.
\end{equation*}
\end{corollary}
\begin{lemma}[Local quasi-boundedness of $N$]\label{lem510} Let $u$ be a non-zero
solution to (\ref{fnls}) with lifespan $I$ that is almost periodic
modulo $G$ with frequency scale function $N: I\rightarrow\Bbb R^+$.
If $K$ is any compact subset of $I$, then
\begin{equation*}
0<\inf_{t\in K}N(t)\leq\sup_{t\in K}N(t)<\infty.
\end{equation*}
\end{lemma}
{\it Proof.} We only prove the first inequality; the argument for
the last is similar.

We argue by contradiction. Suppose that the first inequality fails.
Then, there exists a sequence $t_n\in K$ such that
$\lim_{n\rightarrow\infty}N(t_n)=0$ and hence by Definition
\ref{d1}, $u(t_n)$ converge weakly to zero. Since $K$ is compact, we
can assume $t_n$ converge to a limit $t_0\in K$. As $u\in
C_t^0\dot{H}_x^2(K\times\Bbb R^d)$, we see that $u(t_n)$ converge
strongly to $u(t_0)$. Thus $u(t_0)$ must be zero, contradicting the
hypothesis.
\section{Concentration compactness}
\begin{theorem}[Linear profile decomposition]\label{cc}
Fix $d\geq5$ and $\{u_n\}_{n\geq1}$ be a sequence of functions
bounded in $\dot{H}^2_x(\Bbb R^d)$. Then after passing to a
subsequence if necessary, there exist a sequence of functions
$\{\phi^j\}_{j\geq1}\subset\dot{H}^2_x(\Bbb R^d)$, group elements
$g_n^j\in G$ and times $t_n^j\in \Bbb R$ such that we have the
decomposition
\begin{equation}\label{equ61}
u_n=\sum_{j=1}^Jg_n^je^{it_n^j\Delta^2}\phi^j+w_n^J
\end{equation}
for all $J>1$; here $w_n^J\in\dot{H}^2_x(\Bbb R^d)$ obey
\begin{equation}\label{equ62}
\lim_{J\rightarrow\infty}\limsup_{n\rightarrow\infty}\big\|e^{it\Delta^2}w_n^J\big\|_{L_{t,x}^\frac{2(d+4)}{d-4}(\Bbb
R\times\Bbb R^d)}=0.
\end{equation}
Moreover, for any $j\neq j'$,
\begin{equation}\label{equ63}
\lim_{n\rightarrow\infty}\Big(\frac{\lambda_n}{\lambda_n'}+\frac{\lambda_n'}{\lambda_n}
+\frac{|t_n\lambda_n^4-t_n'(\lambda_n')^4|}{\lambda_n^2\lambda_n'^2}+\frac{|x_n-x_n'|^2}{\lambda_n\lambda_n'}\Big)=+\infty.
\end{equation}
Furthermore, for any $J\geq1$ we have the kinetic energy decoupling
property
\begin{equation}\label{equ64}
\lim_{n\rightarrow\infty}\Big[\|\Delta
u_n\|_2^2-\sum_{j=1}^J\|\Delta \phi^j\|_2^2-\|\Delta
w_n^J\|_2^2\Big]=0.
\end{equation}
\end{theorem}
\begin{remark}
In fact, for any $(q, r)$ ($q\neq2$) such that
$\frac{4}{q}+\frac{d}{r}=\frac{d}{2}-2$, we have by H\"older's
inequality,
\begin{equation}\label{equ66}
\lim_{J\rightarrow\infty}\limsup_{n\rightarrow\infty}\big\|e^{it\Delta^2}w_n^J\big\|_{L_t^qL_x^r(I\times\Bbb
R^d)}=0.
\end{equation}
Moreover, by interpolation we have
\begin{equation}\label{equ67}
\lim_{J\rightarrow\infty}\limsup_{n\rightarrow\infty}\big\|\nabla
e^{it\Delta^2}w_n^J\big\|_{L_t^aL_x^b(I\times\Bbb R^d)}=0,
\end{equation}
where $\frac{4}{a}+\frac{d}{b}=\frac{d}{2}-1$, $a\neq2$.
\end{remark}
The proof of Theorem \ref{cc} is very similar to Theorem 1.6 of
\cite{sh1}. We need only establish the following
\begin{lemma}\label{lem61}Fix $d\geq5$.
For every $f\in \dot{H}_x^2(\Bbb R^d)$, we have
\begin{equation*}
\|f\|_{L^\frac{2d}{d-4}}\leq C\|\Delta
f\|_{L^2}^\frac{d-4}{d}\|\Delta
f\|_{\dot{B}_{2,\infty}^0}^\frac{4}{d}.
\end{equation*}
\end{lemma}
\begin{lemma}\label{lem62} Let $\{t^j\}$, $\{\lambda^j\}$, $\{x^j\}$ be sequences as
in (\ref{equ63}) and $V^j\in L_{t,x}^\frac{2(d+4)}{d-4}(\Bbb
R\times\Bbb R^d)$ for every $j\geq1$, then
\begin{equation}\label{equ65}
\lim_{n\rightarrow\infty}\Big\|\sum_{j=1}^J\frac{1}{(\lambda_n^j)^\frac{d-4}{2}}V^j\Big(\frac{\cdot-t_n^j}{(\lambda_n^j)^4},
\frac{\cdot-x_n^j}{\lambda_n^j}\Big)\Big\|_{L_{t,x}^\frac{2(d+4)}{d-4}(\Bbb
R\times\Bbb R^d)}^\frac{2(d+4)}{d-4}\leq
\sum_{j=1}^J\|V^j\|_{L_{t,x}^\frac{2(d+4)}{d-4}(\Bbb R\times\Bbb
R^d)}^\frac{2(d+4)}{d-4}.
\end{equation}
\end{lemma}
and \begin{lemma}\label{lem63} For all $J\geq1$ and all $1\leq j\leq
J$, the sequence $e^{-it_n^j\Delta^2}[(g_n^j)^{-1}w_n^J]$ converges
weakly to zero in $\dot{H}_x^2(\Bbb R^d)$ as $n\rightarrow\infty$.
In particular, this implies the kinetic energy decoupling
(\ref{equ64}).
 \end{lemma}

{\it Proof of Lemma \ref{lem61}.} This is a direct adaption of the
proof in Bahouri and Gerard \cite{BG}. For every $A>0$, we decompose
$f=P_{>A}f+P_{\leq A}f$, then we have
\begin{equation*}
\|P_{\leq A}f\|_{L^\infty}\leq\sum_{2^k\leq
A}\|P_kf\|_{L^\infty}\lesssim\sum_{2^k\leq
A}{2^{k(\frac{d}{2}-2)}}\|\Delta P_kf\|_{L^2}\lesssim
A^{\frac{d}{2}-2}\|\Delta
f\|_{\dot{B}^0_{2,\infty}}\triangleq\frac{\lambda}{2}.
\end{equation*}
Then $A(\lambda)=\Big(\frac{\lambda}{2C\|\Delta
f\|_{\dot{B}_{2,\infty}^0}}\Big)^\frac{2}{d-4}$ and
\begin{align*}
m\{|f|>\lambda\}\leq&\ m\{|P_{>A(\lambda)}f|>\frac{\lambda}{2}\}\\
\leq&\ \frac{4}{\lambda^2}\|P_{>A(\lambda)}f\|_2^2\\
\leq&\
\frac{4}{\lambda^2}A(\lambda)^{-4}\big\||\cdot|^2\widehat{P_{>A(\lambda)}f}(\cdot)\big\|_2^2.
\end{align*}
Therefore,
\begin{align*}
\|f\|_{L^\frac{2d}{d-4}}^\frac{2d}{d-4}=&\
\frac{2d}{d-4}\int_0^{+\infty}\lambda^\frac{d+4}{d-4}m\{|f|>\lambda\}d\lambda\\
\lesssim&\
\int_0^\infty\frac{4}{\lambda^2}A(\lambda)^{-4}\lambda^\frac{d+4}{d-4}\Big(\int_{|\xi|>A(\lambda)}|\xi|^4|\hat{f}(\xi)|^2d\xi\Big)d\lambda\\
\lesssim&\ \|\Delta f\|_{\dot{B}_{2,\infty}^0}^\frac{8}{d-4}\|\Delta
f\|_{L^2}^2.
\end{align*}
This completes the proof of Lemma \ref{lem61}.

 {\it Proof of Lemma
\ref{lem62}.} Let
$$V_n^j=\frac{1}{(\lambda_n^j)^\frac{d-4}{2}}V^j\Big(\frac{\cdot-t_n^j}{(\lambda_n^j)^4},
\frac{\cdot-x_n^j}{\lambda_n^j}\Big),$$ it suffices to prove that
\begin{equation*}
\lim_{n\rightarrow\infty}\Big\|\sum_{j=1}^JV_n^j\Big\|_{L_{t,x}^\frac{2(d+4)}{d-4}(\Bbb
R\times\Bbb R^d)}^\frac{2(d+4)}{d-4}\leq
\sum_{j=1}^J\|V^j\|_{L_{t,x}^\frac{2(d+4)}{d-4}(\Bbb R\times\Bbb
R^d)}^\frac{2(d+4)}{d-4}.
\end{equation*}
We denote the maximal integer less than $a$ by $[a]$ and let
$k(d)=\big[\frac{2(d+4)}{d-4}\big]$, then
\begin{align*}
\Big\|\sum_{j=1}^JV_n^j&\Big\|_{L_{t,x}^\frac{2(d+4)}{d-4}(\Bbb
R\times\Bbb R^d)}^\frac{2(d+4)}{d-4}=\ \int_{\Bbb R\times\Bbb
R^d}\big|\sum_{j=1}^JV_n^j\big|^{k(d)}\big|\sum_{j=1}^JV_n^j\big|^{\frac{2(d+4)}{d-4}-k(d)}dtdx\\
&\leq\ \sum_{j=1}^J\int_{\Bbb R\times\Bbb
R^d}\big|V_n^j\big|^{k(d)}\big|\sum_{j=1}^JV_n^j\big|^{\frac{2(d+4)}{d-4}-k(d)}dtdx\\
&\ +\sum_{j=1}^J\sum_{j'\neq j}\sum_{j_3,\cdots, j_{k(d)}}\int_{\Bbb
R\times\Bbb R^d}\big|V_n^jV_n^{j'}\big|\big|V_n^{j_3}\cdots
V_n^{J_{k(d)}}\big|\big|\sum_{j=1}^JV_n^j\big|^{\frac{2(d+4)}{d-4}-k(d)}dtdx\\
&=A+B
\end{align*}
We estimate $A$ first.
\begin{equation*}
A\leq\sum_{j=1}^J\int_{\Bbb R\times\Bbb
R^d}\big|V_n^j\big|^\frac{2(d+4)}{d-4}dtdx+\sum_{j=1}^J\sum_{j\neq
j'}\int_{\Bbb R\times\Bbb
R^d}\big|V_n^j\big|^{k(d)}\big|V_n^{j'}\big|^{\frac{2(d+4)}{d-4}-k(d)}dtdx.
\end{equation*}
The second term can be written as
\begin{align*}
&\sum_{j=1}^J\sum_{j\neq j'}\int_{\Bbb R\times\Bbb
R^d}\big|V_n^j\big|^{2k(d)-\frac{2(d+4)}{d-4}}\big|V_n^jV_n^{j'}\big|^{\frac{2(d+4)}{d-4}-k(d)}\\
\lesssim&\ \sum_{j=1}^J\sum_{j\neq
j'}\big\|V_n^j\big\|_{L_{t,x}^\frac{2(d+4)}{d-4}(\Bbb R\times\Bbb
R^d)}^{2k(d)-\frac{2(d+4)}{d-4}}\big\|V_n^jV_n^{j'}\big\|_{L_{t,x}^\frac{d+4}{d-4}(\Bbb
R\times\Bbb R^d)}^{\frac{2(d+4)}{d-4}-k(d)}\\
=&\ \sum_{j=1}^J\sum_{j\neq
j'}\big\|V^j\big\|_{L_{t,x}^\frac{2(d+4)}{d-4}(\Bbb R\times\Bbb
R^d)}^{2k(d)-\frac{2(d+4)}{d-4}}\big\|V_n^jV_n^{j'}\big\|_{L_{t,x}^\frac{d+4}{d-4}(\Bbb
R\times\Bbb R^d)}^{\frac{2(d+4)}{d-4}-k(d)},
\end{align*}
which is $o(1)$ as $n\rightarrow\infty$ by (\ref{equ511}). Now we
estimate $B$. By H\"older's inequality,
\begin{align*}
&\int_{\Bbb R\times\Bbb
R^d}\big|V_n^jV_n^{j'}\big|\big|V_n^{j_3}\cdots
V_n^{J_{k(d)}}\big|\big|\sum_{j=1}^JV_n^j\big|^{\frac{2(d+4)}{d-4}-k(d)}dtdx\\
\leq&\ \sum_{l=1}^J\big\|
V_n^jV_n^{j'}\big\|_{L_{t,x}^\frac{d+4}{d-4}(\Bbb R\times\Bbb
R^d)}\big\|V_n^{j_3}\big\|_{L_{t,x}^\frac{2(d+4)}{d-4}(\Bbb
R\times\Bbb R^d)}
\cdots\big\|V_n^{j_{k(d)}}\big\|_{L_{t,x}^\frac{2(d+4)}{d-4}(\Bbb
R\times\Bbb R^d)}\big\|V_n^l\big\|_{L_{t,x}^\frac{2(d+4)}{d-4}(\Bbb
R\times\Bbb
R^d)}\\
&\longrightarrow0\ \ \text{as}\ \ n\rightarrow\infty.
\end{align*}
Thus we establish Lemma \ref{lem62}.

{\it Proof of Lemma \ref{lem63}.} Fix $J\geq1$ and $1\leq j\leq J$.
Notice that $\{u_n\}_{n\geq1}$ and $\phi^j$ are bounded in
$\dot{H}_x^2(\Bbb R^d)$, by (\ref{equ61}) we deduce that
$\{e^{-it_n^j\Delta^2}[(g_n^j)^{-1}w_n^J]\}_{n\geq1}$ is bounded in
$\dot{H}_x^2(\Bbb R^d)$. Using Alaoglu's Theorem (and passing to a
subsequence if necessary), we obtain that
$e^{-it_n^j\Delta^2}[(g_n^j)^{-1}w_n^J]$ converges weakly in
$\dot{H}_x^2(\Bbb R^d)$ to some $\psi\in \dot{H}_x^2(\Bbb R^d)$. To
prove this lemma, it suffices to show that $\psi\equiv0$.

By weak convergence and (\ref{equ61}),
\begin{align*}
\|\psi\|_{\dot{H}_x^2(\Bbb
R^d)}=&\lim_{n\rightarrow\infty}\langle\Delta
e^{-it_n^j\Delta^2}[(g_n^j)^{-1}w_n^J], \Delta\psi\rangle\\
=&\lim_{n\rightarrow\infty}\langle\Delta
e^{-it_n^j\Delta^2}[(g_n^j)^{-1}(\sum_{l=J+1}^Lg_n^le^{it_n^l\Delta^2}\varphi^l+w_n^L)],\Delta\psi\rangle\\
=&\sum_{l=J+1}^L\lim_{n\rightarrow\infty}\langle\Delta
g_n^le^{it_n^l\Delta^2}\varphi^l,\Delta
g_n^je^{it_n^j\Delta^2}\psi\rangle\\
&\quad+\lim_{n\rightarrow\infty}\langle\Delta
e^{-it_n^j\Delta^2}(g_n^j)^{-1}w_n^L, \Delta\psi\rangle,
\end{align*}
for all $L>J$. By (\ref{equ512}), \begin{equation*}
\lim_{n\rightarrow\infty}\langle\Delta
g_n^le^{it_n^l\Delta^2}\varphi^l,\Delta
g_n^je^{it_n^j\Delta^2}\psi\rangle=0
\end{equation*}
for all $L\geq l\geq J+1>j$.

On the other hand, combining the fact that the family
$\{e^{-it_n^j\Delta^2}[(g_n^j)^{-1}w_n^L]\}_{n, L\geq1}$ is bounded
in $\dot{H}^2_x(\Bbb R^d)$ with
\begin{equation*}
\lim_{L\rightarrow\infty}\limsup_{n\rightarrow\infty}S_{\Bbb
R}\big(e^{it\Delta^2}e^{-it_n^j\Delta^2}[(g_n^j)^{-1}w_n^L]\big)
=\lim_{L\rightarrow\infty}\limsup_{n\rightarrow\infty}S_{\Bbb
R}(e^{it\Delta^2}w_n^L)=0,
\end{equation*}
we deduce that $e^{-it_n^j\Delta^2}(g_n^j)^{-1}w_n^L$ converges
weakly to zero in $\dot{H}_x^2(\Bbb R^d)$ as $n,
L\rightarrow\infty$. Thus for $L$ sufficiently large,
\begin{equation*}
\limsup_{n\rightarrow\infty}\big|\langle\Delta
e^{-it_n^j\Delta^2}(g_n^j)^{-1}w_n^L,
\Delta\psi\rangle\big|\leq\frac{1}{2}\|\psi\|_{\dot{H}_x^2(\Bbb
R^d)}^2.
\end{equation*}
So we have $\psi\equiv0$. This finishes the proof of Lemma
\ref{lem63}.
\section{Reduction to almost periodic
solutions}
\begin{theorem}\label{raps}
Suppose $d\geq5$ is such that Conjecture \ref{cj} failed. Then there
exists a maximal-lifespan solution $u: I\times\Bbb
R^d\rightarrow\Bbb C$ to (\ref{fnls}) such that $\sup_{t\in
I}\|\Delta u(t)\|_2<\|\Delta W\|_2$, $u$ is almost periodic modulo
$G$, and $u$ blows up both forward and backward in time. Moreover,
$u$ has minimal kinetic energy among all blowup solutions, that is
\begin{equation*}
\sup_{t\in I}\|\Delta u(t)\|_2\leq\sup_{t\in J}\|\Delta v(t)\|_2
\end{equation*}
for all maximal-lifespan solution $v: J\times\Bbb R^d\rightarrow\Bbb
C$ that blows up at least one time direction. If furthermore
$d\geq5$ and Conjecture \ref{cj} failed for spherically symmetric
data, then we can also ensure that $u$ is spherically symmetric and
almost periodic modulo $G_{rad}$.
\end{theorem}

For any $0\leq E_0\leq \|\Delta W\|_2^2$, we define
\begin{equation*}
L(E_0):=\sup\{S(u): u:I\times\Bbb R^d\rightarrow\Bbb C\ \text{such
that}\ \sup_{t\in I}\|\Delta u(t)\|_2^2\leq E_0\},
\end{equation*}
where the supremum is taken over all solutions $u:I\times\Bbb
R^d\rightarrow\Bbb C$ to (\ref{fnls}) obeying $\sup_{t\in I}\|\Delta
u(t)\|_2^2\leq E_0$. Thus $L: [0, \|\Delta W\|_2^2]\rightarrow[0,
\infty]$ is a non-decreasing function with $L(\|\Delta
W\|_2^2)=\infty$. Moreover, from Theorem \ref{lw},
\begin{equation*}
L(E_0)\lesssim_d E_0^\frac{d+4}{d-4}\ \ \text{for}\ E_0\leq \eta_0,
\end{equation*}
where $\eta_0=\eta_0(d)$ is the threshold from the small data
theory.

From Theorem \ref{stability}, we see that $L$ is continuous.
Therefore, there must exist a unique critical kinetic energy $E_c$
such that $L(E_0)<\infty$ for $E_0<E_c$ and $L(E_0)=\infty$ for
$E_0\geq E_c$. In particular, if $u: I\times\Bbb R^d\rightarrow\Bbb
C$ is a maximal-lifespan solution to (\ref{fnls}) such that
$\sup_{t\in I}\|\Delta u(t)\|_2^2< E_c$, then $u$ is global and
\begin{equation*}
S(u)\leq L(\sup_{t\in I}\|\Delta u(t)\|_2^2).
\end{equation*}
Failure of Conjecture \ref{cj} is equivalent to the existence of
$0<E_c<\|\Delta W\|_2^2$.
\begin{proposition}[Palais-Smale condition modulo symmetries]\label{ps}
Fix $d\geq5$. Let $u_n: I_n\times\Bbb R^d\rightarrow\Bbb C$ be a
sequence of solutions to (\ref{fnls}) such that
\begin{equation}\label{equ71}
\limsup_{n\rightarrow\infty}\sup_{t\in I_n}\|\Delta u_n(t)\|_2^2=E_c
\end{equation}
and let $t_n\in I_n$ be a sequence of times such that
\begin{equation*}
\lim_{n\rightarrow\infty}S_{\geq
t_n}(u_n)=\lim_{n\rightarrow\infty}S_{\leq t_n}(u_n)=\infty.
\end{equation*}
Then the sequence $u_n(t_n)$ has a subsequence which converges in
$\dot{H}^2_x(\Bbb R^d)$ modulo $G$.
\end{proposition}
{\it Proof.} By the time translation symmetry of (\ref{fnls}), we
may set $t_n=0$ for all $n\geq1$. Thus,
\begin{equation}\label{equ722}
\lim_{n\rightarrow\infty}S_{\geq
0}(u_n)=\lim_{n\rightarrow\infty}S_{\leq 0}(u_n)=\infty.
\end{equation}
Applying Theorem \ref{cc} to the sequence $u_n(0)$ (which is bounded
in $\dot{H}_x^2$ by (\ref{equ71})) and passing to a subsequence if
necessary, we obtain the decomposition
\begin{equation*}
u_n(0)=\sum_{j=1}^Jg_n^je^{it_n^j\Delta^2}\phi^j+w_n^J.
\end{equation*}

Redefining the subsequence once for every $j$ and using a diagonal
argument, we may assume that for each $j$, the sequence
$\{t_n^j\}_{n\geq1}$ converges to some $t^j\in[-\infty, +\infty]$.
If $t^j\in (-\infty, +\infty)$, then by replacing $\phi^j$ by
$e^{it^j\Delta^2}\phi^j$ and $t_n^j-t^j$ by $t^j_n$, we may assume
that $t^j=0$. Moreover, absorbing the error $\sum_{1\leq j\leq J:
t^j=0}g_n^j(e^{it_n^j\Delta^2}\phi^j-\phi^j)$ into the error term
$w_n^J$, we may assume that $t_n^j\equiv0$. Thus either
$t_n^j\equiv0$ or $t_n^j\rightarrow\pm\infty$.

We now define the nonlinear profiles $v^j: I^j\times\Bbb
R^d\rightarrow\Bbb C$ associated to $\phi^j$ and $t_n^j$ as follows:
\renewcommand{\labelenumi}{$\bullet$}
\begin{enumerate}
\item If $t_n^j\equiv0$, then $v^j$ is the maximal-lifespan solution
to (\ref{fnls}) with initial data $v^j(0)=\phi^j$.
\item If $t_n^j\rightarrow+\infty$, then $v^j$ is the
maximal-lifespan solution to (\ref{fnls}) that scatters forward in
time to $e^{it\Delta^2}\phi^j$.
\item If $t_n^j\rightarrow-\infty$, then $v^j$ is the
maximal-lifespan solution to (\ref{fnls}) that scatters backward in
time to $e^{it\Delta^2}\phi^j$.
\end{enumerate}

For each $j,n\geq1$, we define $v_n^j: I_n^j\times\Bbb
R^d\rightarrow\Bbb C$ by
\begin{equation*}
v_n^j(t):=T_{g_n^j}[v^j(\cdot+t_n^j)](t),
\end{equation*}
where $I_n^j:=\{t\in\Bbb R: (\lambda_n^j)^{-4}t+t_n^j\in I^j\}$.
Each $v_n^j$ is a solution to (\ref{fnls}) with initial data at time
$t=0$ given by $v_n^j(0)=g_n^jv^j(t_n^j)$ and maximal lifespan
$I_n^j=(-T_{n,j}^-, T_{n,j}^+)$, where
$-\infty\leq-T_{n,j}^-<0<T_{n,j}^+\leq+\infty$.

By (\ref{equ64}), there exists $J_0\geq1$ such that
\begin{equation*}
\|\Delta \phi^j\|_2\leq\eta_0\ \ \text{for all}\ \ j\geq J_0,
\end{equation*}
where $\eta_0=\eta_0(d)$ is the threshold for the small data theory.
Here, by Theorem \ref{lw} for all $n\geq1$ and all $j\geq J_0$ the
solution $v_n^j$ are global and moreover,
\begin{equation}\label{equ72}
\sup_{t\in\Bbb R}\|\Delta v_n^j\|_2^2+S_{\Bbb R}(v_n^j)\lesssim
\|\Delta \phi^j\|_2^2.
\end{equation}
\begin{lemma}[At least one bad profile]\label{onebad}
There exists $1\leq j_0< J_0$ such that
\begin{equation*}
\limsup_{n\rightarrow\infty}S_{[0, T_{n, j_0}^+)}(v_n^{j_0})=\infty.
\end{equation*}
\end{lemma}

 {\it Proof.} Assume for a contradiction
that for all $1\leq j<J_0$,
\begin{equation}\label{equ73}
\limsup_{n\rightarrow\infty}S_{[0, T_{n, j}^+)}(v_n^j)<\infty,
\end{equation}
which implies $T_{n, j}^+=\infty$ for all $1\leq j<J_0$ and all
sufficiently large $n$. Moreover, subdividing $[0, +\infty)$ into
intervals where the scattering size of $v_n^j$ is small, applying
the Strichartz inequality on each such interval, and then summing,
we obtain
\begin{equation*}
\limsup_{n\rightarrow\infty}\|v_n^j\|_{\dot{S}^2([0,
\infty))}<\infty
\end{equation*}
for all $1\leq j<J_0$, where
$$\|f\|_{\dot{S}^2(I)}:=\sup_{(q,r): \frac{4}{q}+\frac{d}{r}=
\frac{d}{2}}\|\Delta f\|_{L_t^qL_x^r(I\times\Bbb R^d)}.$$ By
(\ref{equ71}), (\ref{equ72}), (\ref{equ73}) and (\ref{equ64}), we
have
\begin{equation}\label{equ74}
\sum_{j\geq1}S_{[0, +\infty)}(v_n^j)\lesssim 1+\sum_{j\geq
J_0}\|\Delta \phi^j\|_2^2\lesssim 1+E_c
\end{equation}
for all $n$ large enough. Now we define the approximation
\begin{equation*}
u_n^J(t):=\sum_{j=1}^Jv_n^j(t)+e^{it\Delta^2}w_n^J.
\end{equation*}
Note that
\begin{align*}
\|u_n^J(0)-u_n(0)\|_{\dot{H}^2_x(\Bbb R^d)}\lesssim&\
\|\sum_{j=1}^J\big(g_n^jv^j(t_n^j)-g_n^je^{it_n^j\Delta^2}\phi^j\big)\|_{\dot{H}^2_x(\Bbb
R^d)}\\
\lesssim&\
\sum_{j=1}^J\|v^j(t_n^j)-e^{it_n^j\Delta^2}\phi^j\|_{\dot{H}^2_x(\Bbb
R^d)}
\end{align*}
and hence, by our choice of $v^j$,
\begin{equation*}
\limsup_{n\rightarrow\infty}\|u_n^J(0)-u_n(0)\|_{\dot{H}^2_x(\Bbb
R^d)}=0.
\end{equation*}
We now show that $u_n^J$ does not blowup forward in time. Indeed, by
(\ref{equ511}), the fact that $v_n^j$ does not blow up forward in
time, Lemma \ref{lem62}, (\ref{equ62}) and (\ref{equ74}), we have
\begin{align}\label{equ75}
\lim_{J\rightarrow\infty}\limsup_{n\rightarrow\infty}S_{[0,\infty)}(u_n^J)\lesssim&\
\lim_{J\rightarrow\infty}\limsup_{n\rightarrow\infty}\Big(S_{[0,\infty)}(\sum_{j=1}^Jv_n^j)+S_{[0,\infty)}(e^{it\Delta^2}w_n^J)\Big)\nonumber\\
\lesssim&\
\lim_{J\rightarrow\infty}\limsup_{n\rightarrow\infty}\sum_{j=1}^JS_{[0,\infty)}(v_n^j)\lesssim\
1+E_c.
\end{align}
Similarly, we can obtain that
\begin{equation}\label{equ76}
\lim_{J\rightarrow\infty}\limsup_{n\rightarrow\infty}\big\|u_n^J\big\|_{\dot{S}^2([0,
\infty))}\leq C(E_c)<\infty.
\end{equation}
In order to apply Theorem \ref{stability}, it suffices to show
$u_n^J$ asymptotically solves (\ref{fnls}) in the sense that
\begin{equation}\label{equ77}
\lim_{J\rightarrow\infty}\limsup_{n\rightarrow\infty}\|\nabla
[(i\partial_t+\Delta^2)u_n^J-F(u_n^J)]\|_{L_t^2L_x^\frac{2d}{d+2}([0,
\infty)\times\Bbb R^d)}=0,
\end{equation}
which reduces to proving
\begin{equation}\label{equ711}
\lim_{J\rightarrow\infty}\limsup_{n\rightarrow\infty}\Big\|\nabla\Big(\sum_{j=1}^JF(v_n^j)-F(\sum_{j=1}^Jv_n^j)\Big)\Big\|_{L_t^2L_x^\frac{2d}{d+2}([0,
\infty)\times\Bbb R^d)}=0
\end{equation}
and
\begin{equation}\label{equ712}
\lim_{J\rightarrow\infty}\limsup_{n\rightarrow\infty}\Big\|\nabla\Big(F(u_n^J-e^{it\Delta^2}w_n^J)-F(u_n^J)\Big)\Big\|_{L_t^2L_x^\frac{2d}{d+2}([0,
\infty)\times\Bbb R^d)}=0.
\end{equation}

We first address (\ref{equ711}). Note that we can write
\begin{equation*}
\Big|\nabla \big[\sum_{j=1}^JF(v_n^j)-F(\sum_{j=1}^J
v_n^j)\big]\Big|\lesssim_J \sum_{j\neq j'}|\nabla
v_n^j||v_n^{j'}|^\frac{8}{d-4}
\end{equation*}
for $d>12$ and
\begin{equation*}
\Big|\nabla \big[\sum_{j=1}^JF(v_n^j)-F(\sum_{j=1}^J
v_n^j)\big]\Big|\lesssim_J \sum_{j\neq j'}\big(|\nabla
v_n^j||v_n^{j'}|^\frac{8}{d-4}+|v_n^j||v_n^{j'}|^\frac{12-d}{d-4}|\nabla
v_n^{j'}|\big)
\end{equation*}
for $5\leq d\leq12$. By the similar argument deriving
(\ref{equ511}), for any $j\neq j'$, we have
\begin{equation*}
\limsup_{n\rightarrow\infty}\big\||v_n^{j'}|^\frac{8}{d-4}|\nabla
v_n^j|\big\|_{L_t^2L_x^\frac{2d}{d+2}([0, \infty)\times\Bbb R^d)}=0
\end{equation*}
for all $d\geq5$ and
\begin{equation*}
\limsup_{n\rightarrow\infty}\big\||v_n^j||v_n^{j'}|^\frac{12-d}{d-4}|\nabla
v_n^{j'}||\big\|_{L_t^2L_x^\frac{2d}{d+2}([0, \infty)\times\Bbb
R^d)}=0
\end{equation*}
for $5\leq d\leq12$. Thus we have
\begin{align*}
&\limsup_{n\rightarrow\infty}\Big\|\nabla
\big[\sum_{j=1}^JF(v_n^j)-F(\sum_{j=1}^J
v_n^j)\big]\Big\|_{L_t^2L_x^\frac{2d}{d+2}([0, \infty)\times\Bbb
R^d)}\\
\lesssim_J& \limsup_{n\rightarrow\infty}\sum_{j\neq
j'}\Big\||v_n^{j'}|^\frac{8}{d-4}|\nabla
v_n^j|\Big\|_{L_t^2L_x^\frac{2d}{d+2}([0, \infty)\times\Bbb R^d)}=0
\end{align*}
and (\ref{equ711}) follows.

We now consider (\ref{equ712}). In dimensions $d\geq12$, by H\"older
and interpolation, we have
\begin{align*}
&\Big\|\nabla\Big(F(u_n^J-e^{it\Delta^2}w_n^J)-F(u_n^J)\Big)\Big\|_{L_t^2L_x^\frac{2d}{d+2}([0,
\infty)\times\Bbb R^d)}\\
\lesssim&\
\big\|e^{it\Delta^2}w_n^J\big\|_{L_t^\frac{8(d+4)}{3(d-4)}L_x^\frac{2d(d+4)}{(d+1)(d-4)}([0,
\infty)\times\Bbb R^d)}^\frac{8}{d-4}\big\|\nabla
u_n^J\big\|_{L_{t,x}^\frac{2(d+4)}{d-2}([0, \infty)\times\Bbb
R^d)}\\
&\
+\big\|e^{it\Delta^2}w_n^J\big\|_{L_t^\frac{8(d+4)}{3(d-4)}L_x^\frac{2d(d+4)}{(d+1)(d-4)}([0,
\infty)\times\Bbb R^d)}^\frac{8}{d-4}\big\|\nabla
e^{it\Delta^2}w_n^J\big\|_{L_{t,x}^\frac{2(d+4)}{d-2}([0,
\infty)\times\Bbb R^d)}\\
&\
+\big\|u_n^J\big\|_{L_t^\frac{8(d+4)}{3(d-4)}L_x^\frac{2d(d+4)}{(d+1)(d-4)}([0,
\infty)\times\Bbb R^d)}^\frac{8}{d-4}\big\|\nabla
e^{it\Delta^2}w_n^J\big\|_{L_{t,x}^\frac{2(d+4)}{d-2}([0,
\infty)\times\Bbb R^d)},
\end{align*}
so (\ref{equ712}) follows from (\ref{equ75}), (\ref{equ66})
,(\ref{equ67}) and the fact that $e^{it\Delta^2}w_n^J$ is bounded in
$\dot{S}^2$. In dimensions $5\leq d<12$, one must add the term
\begin{equation*}
\big\|u_n^J\big\|_{L_t^\infty L_x^\frac{2d}{d-4}([0,
\infty)\times\Bbb R^d)}^\frac{12-d}{d-4}\big\|\nabla
u_n^J\big\|_{L_t^\frac{8}{3}L_x^\frac{2d}{d-5}([0, \infty)\times\Bbb
R^d)}\big\|e^{it\Delta^2}w_n^J\big\|_{L_t^8L_x^\frac{2d}{d-5}([0,
\infty)\times\Bbb R^d)},
\end{equation*}
which is acceptable, too.

We are now in a position to apply Theorem \ref{stability}; invoking
(\ref{equ75}), we conclude that for $n$ sufficiently large,
\begin{equation*}
S_{[0, \infty)}(u_n)\lesssim 1+E_c,
\end{equation*}
this contradicts (\ref{equ722}). This finishes the proof of Lemma
\ref{onebad}.

Let us return to the proof the Proposition \ref{ps} now. Rearranging
the indices, we may assume that there exists $1\leq J_1< J_0$ such
that
\begin{equation*}
\limsup_{n\rightarrow\infty}S_{[0, T_{n,j}^+)}(v_n^j)=\infty\ \
\text{for}\ \ 1\leq j\leq J_1
\end{equation*}
and
\begin{equation}\label{equ723}
\limsup_{n\rightarrow\infty}S_{[0, \infty)}(v_n^j)<\infty\ \
\text{for}\ \ j>J_1.
\end{equation}
Passing to a subsequence in $n$, we can guarantee that $S_{[0,
T_{n,1}^+)}(v_n^1)\rightarrow\infty$.

For each $m,n\geq1$ let us define an integer $j(m,n)\in\{1,\cdots,
J_1\}$ and an interval $K_n^m$ of the form $[0, \tau)$ by
\begin{equation}\label{equ724}
\sup_{1\leq j\leq J_1}S_{K_n^m}(v_n^j)=S_{K_n^m}(v_n^{j(m,n)})=m.
\end{equation}
By the pigeonhole principle, there is a $1\leq j_1\leq J_1$, so that
for infinite many $m$, one has $j(m,n)=j_1$, for infinite many $n$.
Note that the infinite set of $n$ for which this holds may be
m-dependent. By reordering the indices, we may assume that $j_1=1$.
Then by the definition of the critical kinetic energy, we obtain
\begin{equation}\label{equ713}
\limsup_{m\rightarrow\infty}\limsup_{n\rightarrow\infty}\sup_{t\in
K_n^m}\|\nabla v_n^1(t)\|_2^2\geq E_c.
\end{equation}
On the other hand, by virtue of (\ref{equ723}) and (\ref{equ724}),
all $v_n^j$ have finite scattering size on $K_n^m$ for each
$m\geq1$. Thus, by the same argument used in Lemma \ref{onebad}, we
see that for $n$ and $J$ sufficiently large, $u_n^J$ is a good
approximation to $u_n$ on each $K_n^m$. More precisely,
\begin{equation}\label{equ714}
\lim_{J\rightarrow\infty}\limsup_{n\rightarrow\infty}\|u_n^J-u_n\|_{L_t^\infty\dot{H}^2(K_n^m\times\Bbb
R^d)}=0
\end{equation}
for each $m\geq1$.
\begin{lemma}[Kinetic energy decoupling for $u_n^J$]\label{ked}
For all $J\geq1$ and $m\geq1$,
\begin{equation*}
\limsup_{n\rightarrow\infty}\sup_{t\in K_n^m}\Big|\|\Delta
u_n^J\|_2^2-\sum_{j=1}^J\|\Delta v_n^j(t)\|_2^2-\|\Delta
w_n^J\|_2^2\Big|=0.
\end{equation*}
\end{lemma}
{\it Proof.} Fix $J\geq1$ and $m\geq1$. Then for all $t\in K_n^m$,
\begin{align*}
\|\Delta u_n^J\|_2^2=&\ \sum_{j=1}^J\|\Delta v_n^j\|_2^2+\|\Delta
w_n^J\|_2^2+\sum_{j\neq j'}\langle\Delta v_n^j(t), \Delta
v_n^{j'}(t)\rangle\\
&\ +\sum_{j=1}^J\big(\langle\Delta e^{it\Delta^2}w_n^J, \Delta
v_n^j(t)\rangle+\langle\Delta v_n^j(t), \Delta
e^{it\Delta^2}w_n^J\rangle\big).
\end{align*}
It suffices to prove that for all sequences $t_n\in K_n^m$,
\begin{equation}\label{equ720}
\lim_{n\rightarrow\infty}\langle\Delta v_n^j(t_n), \Delta
v_n^{j'}(t_n)\rangle=0
\end{equation}
and
\begin{equation}\label{equ721}
\lim_{n\rightarrow\infty}\langle\Delta e^{it_n\Delta^2}w_n^J, \Delta
v_n^j(t_n)\rangle=0
\end{equation}
for all $1\leq j, j'\leq J$ with $j\neq j'$. We will only
demonstrate the latter, which requires Lemma \ref{lem63}; the former
can be deduced in much the same manner using (\ref{equ63}). By a
change of variables,
\begin{equation}\label{equ719}
\langle\Delta e^{it_n^j\Delta^2}w_n^J, \Delta
v_n^j(t_n)\rangle=\langle\Delta
e^{it_n(\lambda_n^j)^{-4}\Delta^2}\big[(g_n^j)^{-1}w_n^J\big],
\Delta v^j\big(\frac{t_n}{(\lambda_n^j)^4}+t_n^j\big)\rangle.
\end{equation}
As $t_n\in K_n^m\subset[0, T_{n,j}^+)$ for all $1\leq j\leq J_1$, we
have $t_n(\lambda_n^j)^{-4}+t_n^j\in I^j$ for all $j\geq1$. Recall
that $I^j$ is the maximal lifespan of $v^j$; for $j>J_1$ we have
$\Bbb R^+\subset I^j$. By refining the sequence once for every $j$
and using the standard diagonalisation argument, we may assume
$t_n(\lambda_n^j)^{-4}+t_n^j$ converges for every $j$.

Fix $1\leq j\leq J$. If $t_n(\lambda_n^j)^{-4}+t_n^j$ converges to
some point $\tau^j$ in the interior of $I^j$, then by the continuity
of the flow, $v^j(t_n(\lambda_n^j)^{-4}+t_n^j)$ converges to
$v^j(\tau^j)$ in $\dot{H}^2_x(\Bbb R^d)$. On the other hand, by
(\ref{equ64}),
\begin{equation*}
\limsup_{n\rightarrow\infty}\big\|e^{it_n(\lambda_n^j)^{-4}\Delta^2}\big[(g_n^j)^{-1}w_n^J\big]\big\|_{\dot{H}^2_x(\Bbb
R^d)}=\limsup_{n\rightarrow\infty}\big\|w_n^J\big\|_{\dot{H}^2_x(\Bbb
R^d)}\lesssim E_c^{\frac{1}{2}}.
\end{equation*}
Combining this with (\ref{equ719}), we obtain
\begin{align*}
\lim_{n\rightarrow\infty}\langle\Delta e^{it_n^j\Delta^2}w_n^J,
\Delta v_n^j(t_n)\rangle=&\lim_{n\rightarrow\infty}\langle\Delta
e^{it_n(\lambda_n^j)^{-4}\Delta^2}\big[(g_n^j)^{-1}w_n^J\big],
\Delta v^j\big(\tau^j\big)\rangle\\
=&\lim_{n\rightarrow\infty}\langle\Delta
e^{-it_n^j\Delta^2}\big[(g_n^j)^{-1}w_n^J\big], \Delta
e^{-i\tau^j\Delta^2}v^j(\tau^j)\rangle.
\end{align*}
Invoking Lemma \ref{lem63}, we deduce (\ref{equ721}).

Consider the case when $t_n(\lambda_n^j)^{-4}+t_n^j$ converges to
$\sup{I^j}$. Then we must have $\sup{I^j}=\infty$ and $v^j$ scatters
forward in time. In fact, this is clearly true if
$t_n^j\rightarrow\infty$ as $n\rightarrow\infty$; in other cases,
failure would imply
\begin{equation*}
\limsup_{n\rightarrow\infty}S_{[0,
t_n]}(v_n^j)=\limsup_{n\rightarrow\infty}S_{[t_n^j,
t_n(\lambda_n^j)^{-4}+t_n^j]}(v^j)=\infty,
\end{equation*}
which contradicts $t_n\in K_n^m$. Therefore, there exists
$\psi^j\in\dot{H}_x^2(\Bbb R^d)$ such that
\begin{equation*}
\lim_{n\rightarrow\infty}\big\|v^j(t_n(\lambda_n^j)^{-4}+t_n^j)-e^{i(t_n(\lambda_n^j)^{-4}+t_n^j)\Delta^2}\psi^j\big\|_{\dot{H}^2_x(\Bbb
R^d)}=0.
\end{equation*}
Together with (\ref{equ719}), this yields
\begin{equation*}
\lim_{n\rightarrow\infty}\langle\Delta e^{it_n^j\Delta^2}w_n^J,
\Delta v_n^j(t_n)\rangle=\lim_{n\rightarrow\infty}\langle\Delta
e^{-it_n^j\Delta^2}\big[(g_n^j)^{-1}w_n^J\big], \Delta
\psi^j\rangle,
\end{equation*}
which by Lemma \ref{lem63} implies (\ref{equ721}).

Finally, we consider the case when $t_n(\lambda_n^j)^{-4}+t_n^j$
converges to $\inf{I^j}$. Since $t_n(\lambda_n^j)^{-4}\geq0$ and
$\inf{I^j}<\infty$ for all $j\geq1$ we see that $t_n^j$ does not
converge to $+\infty$. Moreover, if $t_n^j\equiv0$, then
$\inf{I^j}<0$; as $t_n(\lambda_n^j)^{-4}\geq0$, we see that $t_n^j$
cannot be identically zero. This leaves $t_n^j\rightarrow-\infty$ as
$n\rightarrow\infty$. Thus $\inf{I^j}=-\infty$ and $v^j$ scatters
backward in time to $e^{it\Delta^2}\phi^j$. We obtain
\begin{equation*}
\lim_{n\rightarrow\infty}\big\|v^j(t_n(\lambda_n^j)^{-4}+t_n^j)-e^{i(t_n(\lambda_n^j)^{-4}+t_n^j)\Delta^2}\phi^j\big\|_{\dot{H}^2_x(\Bbb
R^d)}=0,
\end{equation*}
which by (\ref{equ721}) implies
\begin{equation*}
\lim_{n\rightarrow\infty}\langle\Delta e^{it_n^j\Delta^2}w_n^J,
\Delta v_n^j(t_n)\rangle=\lim_{n\rightarrow\infty}\langle\Delta
e^{-it_n^j\Delta^2}\big[(g_n^j)^{-1}w_n^J\big], \Delta
\phi^j\rangle.
\end{equation*}
Invoking Lemma \ref{lem63} once again, we derive (\ref{equ721}).
This finishes the proof of Lemma \ref{ked}.

 Thus by (\ref{equ71})
, (\ref{equ714}) and Lemma \ref{ked}, we have
\begin{equation*}
E_c\geq\limsup_{n\rightarrow\infty}\sup_{t\in K_n^m}\|\Delta
u_n(t)\|_2^2=\lim_{J\rightarrow\infty}\limsup_{n\rightarrow\infty}\big\{\|\Delta
w_n^J\|_2^2+\sup_{t\in K_n^m}\sum_{j=1}^J\|\Delta
v_n^j(t)\|_2^2\big\}.
\end{equation*}
Invoking (\ref{equ713}), this implies $J=1$, $v_n^j\equiv0$ for all
$j\geq2$ and $w_n:=w_n^1$ converges to zero strongly in
$\dot{H}_x^2$. In other words,
\begin{equation}\label{equ715}
u_n(0)=g_ne^{i\tau_n\Delta^2}\phi+w_n
\end{equation}
for some $g_n\in G$, $\tau_n\in\Bbb R$ and some functions $\phi$,
$w_n\in\dot{H}^2_x(\Bbb R^d)$ with $w_n\rightarrow0$ strongly in
$\dot{H}^2_x(\Bbb R^d)$. Moreover, the sequence $\tau_n\equiv0$ or
$\tau_n\rightarrow\pm\infty$.

If $\tau_n\equiv0$, (\ref{equ715}) immediately implies that $u_n(0)$
converges modulo $G$ to $\phi$, which proves Proposition \ref{ps} in
this case.

Finally, we will show that this is the only possible case, that is,
$\tau_n$ cannot converge to either $\infty$ or $-\infty$. We argue
by contradiction. Assume that $\tau_n$ converges to $+\infty$, the
proof in the negative time direction is essentially the same. By the
Strichartz inequality, $S_{\Bbb R}(e^{it\Delta^2}\phi)<\infty$. Thus
we have
\begin{equation*}
\lim_{n\rightarrow\infty}S_{\geq0}(e^{it\Delta^2}e^{i\tau_n\Delta^2}\phi)=0.
\end{equation*}
Since the action of $G$ preserves linear solutions and the
scattering size, this implies
\begin{equation*}
\lim_{n\rightarrow\infty}S_{\geq0}(e^{it\Delta^2}g_ne^{i\tau_n\Delta^2}\phi)=0.
\end{equation*}
Combining this with (\ref{equ715}) and $w_n\rightarrow0$ in
$\dot{H}_x^2$, we conclude
\begin{equation*}
\lim_{n\rightarrow\infty}S_{\geq0}(e^{it\Delta^2}u_n(0))=0.
\end{equation*}
An application of Lemma \ref{stability} yields
\begin{equation*}
\lim_{n\rightarrow\infty}S_{\geq0}(u_n)=0,
\end{equation*}
which contradicts (\ref{equ71}).

{\it Proof of Theorem \ref{raps}.} Suppose $d\geq5$ is such that
Conjecture \ref{cj} failed. Then the critical kinetic energy $E_c$
must obey $E_c<\|\Delta W\|_2^2$. By the definition of the critical
kinetic energy, we can find a sequence $u_n: I_n\times\Bbb
R^d\rightarrow\Bbb C$ of solutions to (\ref{fnls}) with $I_n$
compact,
\begin{equation}\label{equ716}
\sup_{n\geq1}\sup_{t\in I_n}\|\Delta u_n(t)\|_2^2=E_c\ \ \text{and}\
\ \lim_{n\rightarrow\infty}S_{I_n}(u_n)=\infty.
\end{equation}
Let $t_n\in I_n$ be such that $S_{\geq t_n}(u_n)=S_{\leq t_n}(u_n)$.
Then
\begin{equation}\label{equ717}
\lim_{n\rightarrow\infty}S_{\geq
t_n}(u_n)=\lim_{n\rightarrow\infty}S_{\leq t_n}(u_n)=\infty.
\end{equation}
Using the time translation, we may take all $t_n=0$.

Applying Proposition \ref{ps} and passing to a subsequence if
necessary, we can find $g_n\in G$ and a function
$u_0\in\dot{H}^2_x(\Bbb R^d)$ such that $g_nu_n(0)\rightarrow u_0$
strongly in $\dot{H}^2_x(\Bbb R^d)$. By applying the group action
$T_{g_n}$ to the solution $u_n$ we may take all the $g_n$ to be
identity. Thus $u_n(0)$ converges strongly to $u_0$ in
$\dot{H}^2_x(\Bbb R^d)$.

Let $u: I\times\Bbb R^d\rightarrow\Bbb C$ be the maximal-lifespan
solution to (\ref{fnls}) with initial data $u(0)=u_0$. As
$u_n(0)\rightarrow u_0$ in $\dot{H}^2_x(\Bbb R^d)$, Theorem
\ref{stability} shows that $I\subseteq \liminf I_n$ and
\begin{equation*}
\lim_{n\rightarrow\infty}\|u_n-u\|_{L_t^\infty\dot{H}^2_x(K\times\Bbb
R^d)}=0\ \ \text{for all compact}\ K\subset I.
\end{equation*}
Thus by (\ref{equ716}),
\begin{equation}\label{equ718}
\sup_{t\in I}\|\Delta u(t)\|_2^2\leq E_c.
\end{equation}

Next we prove that $u$ blows up both forward and backward in time.
Indeed, if $u$ does not blow up forward in time, then $[0,
\infty)\subset I$ and $S_{\geq0}(u)<\infty$. By Theorem
\ref{stability}, this implies $S_{\geq0}(u_n)<\infty$ for
sufficiently large $n$, which contradicts (\ref{equ717}). A similar
argument proves that $u$ blows up backward in time.

Therefore, by our definition of $E_c$, $\sup_{t\in I}\|\Delta
u(t)\|_2^2\geq E_c$. Combining this with (\ref{equ718}), we obtain
\begin{equation*}
\sup_{t\in I}\|\Delta u(t)\|_2^2=E_c.
\end{equation*}
It remains to show that $u$ is almost periodic modulo $G$. Consider
an arbitrary sequence $\tau_n\in I$. As $u$ blows up in both time
directions
\begin{equation*}
S_{\geq \tau_n}(u)=S_{\leq\tau_n}(u)=\infty.
\end{equation*}
Applying Proposition \ref{ps}, we conclude that $u(\tau_n)$ admits a
convergent subsequence in $\dot{H}^2_x(\Bbb R^d)$ modulo $G$. Thus
the orbit $\{Gu(t): t\in I\}$ is precompact in
$G\backslash\dot{H}^2_x(\Bbb R^d)$.

A direct analogue of Theorem 7.3 in \cite{TVZ} shows that if $u_n$
is a sequence of bounded radial functions in $\dot{H}^2(\Bbb R^d)$,
then there exists a family of radial functions $\varphi^j$,
$j=1,2,\cdots$ in $\dot{H}^2(\Bbb R^d)$ and group elements
$G_n^{(j)}\in G_{rad}'$ for $j,n=1,2,\cdots$ such that we have
decomposition (\ref{equ61}) for all $l=1,2,\cdots$, where $w_n^l\in
\dot{H}^2$ is radial and obeys (\ref{equ62}). Moreover $g_n^j$,
$g_n^{(j')}$ are asymptotically orthogonal in the sense of
(\ref{equ63}) for any $j\neq j'$ and for $l\geq1$, we have the
energy decoupling property (\ref{equ64}). This concludes the proof
of Theorem \ref{raps}.
\section{Three enemies}
\begin{theorem}[Three special scenarios for blowup]\label{enemies}
Fix $d\geq5$ and suppose that Conjecture \ref{cj} fails for this
choice of $d$. Then there exists a minimal kinetic energy,
maximal-lifespan solution $u: I\times\Bbb R^d\rightarrow\Bbb C$,
which is almost periodic modulo $G$, $S_I(u)=\infty$, and obeys
$sup_{t\in I}\|\Delta u\|_2<\|\Delta W\|_2$. If furthermore $d\geq5$
and Conjecture \ref{cj} failed for spherically symmetric data, then
$u$ may be chosen to be spherically symmetric and almost periodic
modulo $G_{rad}$.

With or without spherical symmetry, we can also ensure that the
lifespan $I$ and the frequency scale function $N: I\rightarrow\Bbb
R^+$ match one of the following three scenarios:
\renewcommand{\labelenumi}{\Roman{enumi}.}
\begin{enumerate}
\item {\rm( Finite time blowup.)} We have that either $|\inf
I|<\infty$ or $\sup I<\infty$.
\item {\rm( Soliton-like solution.)} We have $I=\Bbb R$ and $$N(t)=1,\ \  \text{for all}\ \ t\in\Bbb R.$$
\item {\rm(Low-to-high frequency cascade.)} We have $I=\Bbb R$ and $$\inf_{t\in\Bbb R}N(t)\geq1,\ \
\text{and}\ \ \limsup_{t\rightarrow\infty}N(t)=\infty.$$
\end{enumerate}
\end{theorem}
{\it Proof.} The proof is a straightforward adaptation of the
similar proof in Killip, Tao and Visan \cite{KTV2}, \cite{KM}. Let
$v: J\times\Bbb R^d\rightarrow\Bbb C$ denote a minimal kinetic
energy blowup solution whose existence is guaranteed by Theorem
\ref{raps}. We denote the frequency scale function of $v$ by
$N_v(t)$ and spatial center function of $v$ by $x_v(t)$. For any
$T\geq0$, define the quantity
\begin{equation}\label{osc}
\text{osc}(T)=\inf_{t_0\in J}\frac{\sup\{N_v(t): t\in J\ \
\text{and}\ \ |t-t_0|\leq TN_v(t_0)^{-4}\}}{\inf\{N_v(t): t\in J\ \
\text{and}\ \ |t-t_0|\leq TN_v(t_0)^{-4}\}}.
\end{equation}
{\bf Case I:} $\displaystyle\lim_{T\rightarrow\infty}\text{osc}(T)<\infty$.\\

In this case, we can find a finite number $A=A_v$, a sequence $t_n$
of times in $J$, and a sequence $T_n\rightarrow\infty$ such that
\begin{equation*}
\frac{\sup\{N_v(t): t\in J\ \ \text{and}\ \ |t-t_n|\leq
T_nN_v(t_n)^{-4}\}}{\inf\{N_v(t): t\in J\ \ \text{and}\ \
|t-t_n|\leq T_nN_v(t_n)^{-4}\}}<A
\end{equation*}
for all $n$. Note that this, together with Corollary \ref{c58},
implies that
\begin{equation*}
[t_n-T_n/N_v(t_n)^4, t_n+T_n/N_v(t_n)^4]\subset J
\end{equation*}
and
\begin{equation*}
N_v(t)\sim_vN_v(t_n)
\end{equation*}
for all $t$ in this interval.

Let $v^{[t_n]}$ be the normalisation of $v$ at times $t_n$ as in
(\ref{norm}), then $v^{[t_n]}$ is a maximal-lifespan normalized
solution with lifespan
\begin{equation*}
J_n:=\{s\in\Bbb R: t_n+\frac{1}{N_v(t_n)^4}s\in J\}\supset[-T_n,
T_n]
\end{equation*}
and energy $E(v)$. It is almost periodic modulo $G$ with frequency
scale function
\begin{equation*}
N_{v^{[t_n]}}(s):=\frac{1}{N_v(t_n)}N_v\big(t_n+\frac{1}{N_v(t_n)^4}s\big)
\end{equation*}
and compactness modulus function $C$. In particular, we see that
\begin{equation}\label{eque1}
N_{v^{[t_n]}}(s)\sim_v1
\end{equation}
for all $s\in[-T_n, T_n]$.

We now apply Lemma \ref{lem57} and conclude (passing to a
subsequence if necessary) that $v^{[t_n]}$ converge locally
uniformly to a maximal-lifespan solution $u$ with energy $E(v)$
defined on an open interval $I$ containing 0 and which is almost
periodic modulo $G$. As $T_n\rightarrow\infty$, Lemma \ref{lem56}
and (\ref{eque1}) imply that the frequency scale function $N:
I\rightarrow\Bbb R^+$ of $v$ satisfies
\begin{equation*}
0<\inf_{t\in I} N(t)\leq\sup_{t\in I} N(t)<\infty.
\end{equation*}
In particular, by Corollary \ref{c59}, $I$ has no finite endpoints
and hence $I=\Bbb R$. By modifying $C$ by a bounded amount we may
now normalize $N(t)\equiv1$. Thus we have constructed a soliton-like
solution in the sense of Theorem \ref{enemies}.

When $\text{osc}(T)$ is unbounded, we must seek a solution belonging
to one of the remaining two scenarios. We introduce the quantity
\begin{equation*}
a(t_0)=\frac{N_v(t_0)}{\sup\{N_v(t): t\leq
t_0\}}+\frac{N_v(t_0)}{\sup\{N_v(t): t\geq t_0\}}.
\end{equation*}

{\bf Case II:}
$\displaystyle\lim_{T\rightarrow\infty}\text{osc}(T)=\infty$ and
$\displaystyle\inf_{t_0\in J}a(t_0)=0$.

As $\inf_{t_0\in J}a(t_0)=0$, there exists a sequence of times
$t_n\in J$ such that $a(t_n)\rightarrow0$ as $n\rightarrow\infty$.
By the definition of $a$, we can also find times $t_n^-<t_n<t_n^+$
with $t_n^-$, $t_n^+\in J$ such that
\begin{equation*}
\frac{N_v(t_n^-)}{N_v(t_n)}\rightarrow+\infty\ \ \text{and}\ \
\frac{N_v(t_n^+)}{N_v(t_n)}\rightarrow+\infty.
\end{equation*}
Next we choose times $t_n'\in(t_n^-, t_n^+)$ so that
\begin{equation*}
N_v(t_n')\leq 2\inf\{N(t): t\in[t_n^-, t_n^+]\}.
\end{equation*}
In particular, \begin{equation*} N_v(t_n')\sim\inf_{t_n^-\leq t\leq
t_n^+}N_v(t),
\end{equation*}
which allows us to deduce that
\begin{equation*}
\frac{N_v(t_n^-)}{N_v(t_n')}\rightarrow+\infty\ \ \text{and}\ \
\frac{N_v(t_n^+)}{N_v(t_n')}\rightarrow+\infty.
\end{equation*}
We define the rescaled and translated times $s_n^-<0<s_n^+$ by
\begin{equation*}
s_n^\pm:=N_v(t_n')^4(t_n^\pm-t_n')
\end{equation*}
and the normalisations $v^{[t_n']}$ at times $t_n'$ by (\ref{norm}).
These are  normalised maximal-lifespan solutions with lifespans
containing $[s_n^-, s_n^+]$, which are almost periodic modulo $G$
with frequency scale functions
\begin{equation}\label{eque2}
N_{v^{[t_n']}}(s):=\frac{1}{N_v(t_n')}N_v(t_n'+\frac{1}{N_v(t_n')^4}s).
\end{equation}
By the way we choose $t_n'$, we see that
\begin{equation}\label{eque3}
N_{v^{[t_n']}}(s)\gtrsim1
\end{equation}
for all $s_n^-\leq s\leq s_n^+$. Moreover,
\begin{equation}\label{eque4}
N_{v^{[t_n']}}(s_n^\pm)\rightarrow\infty\ \ \text{as}\ \
n\rightarrow\infty
\end{equation}
for either choice of sign.

We now apply Lemma \ref{lem57} and conclude (passing to subsequence
if necessary) that $v^{[t_n']}$ converge locally uniformly to a
maximal-lifespan solution $u$ of energy $E(v)$ defined on an open
interval $I$ containing $0$, which is almost periodic modulo $G$.

Let $N$ be a frequency scale function for $u$. From Lemma
\ref{lem510} we see that $N(t)$ is bounded from above on any compact
set $K\subset I$. From this, Lemma \ref{lem55} and Lemma
\ref{lem56}, we see that $N_{v^{[t_n']}}(t)$ is also bounded from
above, uniformly in $t\in K$, for all sufficiently large $n$
(depending on $K$). As a consequence of this and (\ref{eque4}), we
see that $s_n^-$ and $s_n^+$ cannot be any limit points in $K$; thus
$K\subset [s_n^-, s_n^+]$ for all sufficiently large $n$. Therefore
$s_n^\pm$ converge to the endpoints of $I$. If $\sup(I)<+\infty$ or
$|\inf(I)|<+\infty$, then $u$ blows up in finite time. Otherwise,
$I=\Bbb R$. In this case, we need to show that
\begin{equation*}
\limsup_{t\rightarrow-\infty}N(t)=\limsup_{t\rightarrow+\infty}N(t)=\infty.
\end{equation*}
By time reversal symmetry, it suffices to establish that
$\lim_{t\rightarrow+\infty}N(t)=\infty$. By (\ref{eque3}) and Lemma
\ref{lem56}, we conclude that
\begin{equation*}
\inf_{t\in \Bbb R}N(t)\gtrsim1.
\end{equation*}
Suppose $\lim_{t\rightarrow+\infty}N(t)<\infty$, then $N(t)\sim_u 1$
for all $t\geq0$. We conclude from Lemma \ref{lem56} that for every
$m\geq1$, there exists an $n_m$ such that
\begin{equation*}
N_{v^{[t_{n_m}']}}(t)\sim_u 1
\end{equation*}
for all $0\leq t\leq m$. But by (\ref{osc}) and (\ref{eque2}) this
implies $\text{osc}(\frac{\varepsilon m}{2})\lesssim1$ for all $m$
and some $\varepsilon=\varepsilon(u)>0$ independent of $m$. Note
that $\varepsilon$ is chosen as a lower bound on the quantities
$N(t_{n_m}'')^4/N(t_{n_m}')^4$ where
$t_{n_m}''=t_{n_m}'+\frac{m}{2}N(t_{n_m}')^{-4}$. This contradicts
the hypothesis $\lim_{T\rightarrow\infty}\text{osc}(T)=\infty$ and
so settles Case II.

 {\bf Case III:}
$\displaystyle\lim_{T\rightarrow\infty}\text{osc}(T)=\infty$ and
$\displaystyle\inf_{t_0\in J}a(t_0)>0$.

Let $\varepsilon=\varepsilon(v)>0$ be such that $\inf_{t_0\in
J}a(t_0)\geq 2\varepsilon$. We call a time $t_0$ future-spreading if
$N(t)\leq\varepsilon^{-1}N(t_0)$ for all $t\geq t_0$; we call a time
$t_0$ past-spreading if $N(t)\leq\varepsilon^{-1}N(t_0)$ for all
$t\leq t_0$. Note that every $t_0\in J$ is future-spreading,
past-spreading or possibly both.

We will show that either all sufficiently late times are
future-spreading or that all sufficiently early times are
past-spreading. We only show the first half because the other half
is similar. If this were false, there would be a future-spreading
time $t_0$ and a sequence of past-spreading times $t_n$ that
converges to $\sup(J)$. For sufficiently large $n$, we have $t_n\geq
t_0$. Since $N_v(t_0)\leq\varepsilon^{-1}N_v(t_n)$ and
$N_v(t_n)\leq\varepsilon^{-1}N_v(t_0)$ we see that
\begin{equation*}
N_v(t_n)\sim_vN_v(t_0)
\end{equation*}
for all such $n$. For any $t_0<t<t_n$, we know that $t$ is either
past-spreading or future-spreading; thus we have either
$N_v(t_0)\leq\varepsilon^{-1}N_v(t)$ or
$N_v(t_n)\leq\varepsilon^{-1}N_v(t)$. Also, since $t_0$  is
future-spreading $N_v(t)\leq\varepsilon^{-1}N_v(t_0)$ and $t_n$ is
past-spreading, $N_v(t)\leq\varepsilon^{-1}N_v(t_n)$, we conclude
that
\begin{equation*}
N_v(t)\sim_vN_v(t_0)
\end{equation*}
for all $t_0<t<t_n$; since $t_n$ converges to $\sup(J)$, this claim
in fact holds for all $t_0<t<\sup(J)$. From Corollary \ref{c59} we
see that $v$ does not blow up forward in finite time, that is,
$\sup(J)=\infty$. This implies that
$\lim_{T\rightarrow\infty}\text{osc}(T)<\infty$, a contradiction. We
may now assume that future-spreading occurs for all sufficiently
late times; more precisely, we can find $t_0\in J$ such that all
times $t\geq t_0$ are future-spreading.

Choose $T$ so that $\text{osc}(T)>2\varepsilon^{-1}$. We will now
recursively construct an increasing sequence of times
$\{t_n\}_{n=1}^\infty$ so that
\begin{equation*}
0\leq t_{n+1}-t_n\leq 8TN_v(t_n)^{-4}\ \ \text{and}\ \
N_v(t_{n+1})\leq\frac{1}{2}N_v(t_n).
\end{equation*}
Given $t_n$, set $t_n':=t_n+16TN_v(t_n)^{-4}$. If $N_v(t_n')\leq
\frac{1}{2}N_v(t_n)$ we choose $t_{n+1}=t_n'$ and the properties set
out above follow immediately. If $N_v(t_n')>\frac{1}{2}N_v(t_n)$,
then
\begin{equation}\label{eque5}
J_n:=[t_n'-TN_v(t_n')^{-4}, t_n'+TN_v(t_n')^{-4}]\subseteq[t_n,
t_n+8TN_v(t_n)^{-4}].
\end{equation}
As $t_n$ is future-spreading, this allows us to conclude that
$N_v(t)\leq\varepsilon^{-1}N_v(t_n)$ on $J_n$, but then by the way
$T$ is chosen, we may find $t_{n+1}\in J_n$ so that
$N_v(t_{n+1})\leq N_v(t_n)$. Having obtained a sequence of times
obeying (\ref{eque5}), we may conclude that any subsequential limit
$u$ of $v^{[t_n]}$ is a finite-time blowup solution. To elaborate,
set $s_n=(t_0-t_n)N_v(t_n)^4$ and note that $N_{v^{[t_n]}}(s_n)\geq
2^n$. However $s_n$ is a bounded sequence; indeed,
\begin{equation*}
|s_n|=N(t_n)^4\sum_{k=0}^{n-1}[t_{k+1}-t_k]\leq
8T\sum_{k=0}^{n-1}\frac{N(t_n)^4}{N(t_k)^4}\leq
8T\sum_{k=0}^{n-1}2^{-(n-k)}\leq 8T.
\end{equation*}
In this way, we see that the solution $u$ must blow up at some time
$-8T\leq t<0$.

 This completes the proof of Theorem \ref{enemies}.
\section{Kill the enemies}
\begin{theorem}[No finite-time blowup]\label{kill}
Let $d\geq5$. Then there are no maximal-lifespan radial solutions
$u: I\times\Bbb R^d\rightarrow\Bbb C$ to (\ref{fnls}) that are
almost periodic modulo $G_{rad}$, obey
\begin{equation}\label{equ81}S_I(u)=\infty,
\end{equation}
\begin{equation*} \sup_{t\in I}\|\Delta u(t)\|_{L^2}<\|\Delta
W\|_{L^2}
\end{equation*}
and are such that either $|\inf I|<\infty$ or $\sup I<\infty$.
\end{theorem}
{\it Proof.} Suppose for a contradiction that there existed such a
solution $u$. Without loss of generality, we may assume that $\sup
I<\infty$. Then by Corollary \ref{c59},
\begin{equation}\label{equ82}
\liminf_{t\nearrow\sup I}N(t)=\infty.
\end{equation}
We now show that
\begin{equation}\label{equ83}
\limsup_{t\nearrow\sup I}\int_{|x|\leq R}|u(t,x)|^2dx=0\ \ \text{for
all}\ \ R>0.
\end{equation}
In fact, let $u(t,x)=N(t)^\frac{d-4}{2}v(N(t)x,t)$, then
\begin{align*}
\int_{|x|\leq R}|u(t,x)|^2dx=&\ N(t)^{-4}\int_{|x|\leq
RN(t)}|v(x,t)|^2dx\\
=&\ N(t)^{-4}\int_{|x|\leq \epsilon
RN(t)}|v(x,t)|^2dx+N(t)^{-4}\int_{|x|\leq RN(t)\atop |x|>\epsilon
RN(t)}|v(x,t)|^2dx\\
\doteq&\ A+B.
\end{align*}
By H\"older's inequality, we have
\begin{equation*}
A\leq (N(t))^{-4}(\epsilon RN(t))^4\|v\|_{L^{2^{\#}}}^2\leq
(\epsilon R)^4\|\Delta W\|_2^2.
\end{equation*}
$A$ can be acceptable if we take $\epsilon$ to be sufficiently
small. By H\"older, (\ref{equ82}) and the fact that $u$ is almost
periodic, we have
\begin{equation*}
B\leq R^4\|v\|_{L^{2^{\#}}(|x|>\epsilon RN(t))}^2\rightarrow 0, \ \
\text{as}\ \ t\rightarrow\sup I.
\end{equation*}
Thus (\ref{equ83}) is proved.

For $t\in I$, define
\begin{equation*}
M_R(t):=\int_{\Bbb R^d}\phi\big(\frac{|x|}{R}\big)|u(t,x)|^2dx,
\end{equation*}
where $\phi$ is a smooth, radial function such that $\phi(r)=1$ for
$r\leq1$ and $\phi(r)=0$ for $r\geq2$. By (\ref{equ83}),
\begin{equation}\label{equ84}
\limsup_{t\nearrow\sup I} M_R(t)=0\ \ \text{for all}\ \ R>0.
\end{equation}
On the other hand,
\begin{equation*}
\partial_tM_R(t)=\ -2{\rm
Im}\int\Delta\Big(\phi\big(\frac{|x|}{R}\big)\Big)\bar{u}\Delta
udx-2{\rm
Im}\int\nabla\Big(\phi\big(\frac{|x|}{R}\big)\Big)\cdot\nabla\bar{u}\Delta
udx.
\end{equation*}
So by H\"older and Hardy's inequality, we have
\begin{align*}
|\partial_tM_R(t)|\lesssim&\ \int_{|x|\sim R}\frac{|u||\Delta
u|}{R^2}dx+\int_{|x|\sim R}\frac{|\nabla u||\Delta u|}{R}\\
\lesssim&\ \big\|\frac{u}{|x|^2}\big\|_2\|\Delta
u\|_2+\big\|\frac{|\nabla u|}{|x|}\big\|_2\|\Delta u\|_2\\
\lesssim&\ \|\Delta u\|_2^2\lesssim\ \|\Delta W\|_2^2.
\end{align*}
Thus,
\begin{equation*}
M_R(t_1)=M_R(t_2)+\int_{t_2}^{t_1}\partial_tM_R(t)dt\lesssim
M_R(t_2)+|t_1-t_2|\|\Delta W\|_2^2
\end{equation*}
for all $t_1, t_2\in I$ and $R>0$. Let $t_2\nearrow\sup I$ and
invoking (\ref{equ84}), we have
\begin{equation*}
M_R(t_1)\lesssim |\sup I-t_1|\|\Delta W\|_2^2.
\end{equation*}
Now letting $R\rightarrow\infty$ and using the conservation of mass,
we obtain $u_0\in L_x^2(\Bbb R^d)$. Finally, letting
$t_1\nearrow\sup I$, we deduce $u_0=0$. Thus $u\equiv0$,
contradicting (\ref{equ81}).
\begin{theorem}[Absence of cascades and solitons] Let $d\geq5$.
There are no global radial solutions to (\ref{fnls}) that are
low-to-high cascades or solitons in the sense of Theorem
\ref{enemies}.
\end{theorem}
{\it Proof.} We will show that no global radial solutions that are
almost periodic modulo $G_{rad}$ with the frequency scale function
$N(t)\geq1$ for all $t\in \Bbb R$.

By the almost periodicity of $u$ and Hardy's inequality, we have
that for any $\epsilon>0$, there exists $R(\epsilon)>0$ such that
for all $t\in [0, +\infty)$,
\begin{equation}\label{equ85}
\int_{|x|>R(\epsilon)}\Big(|\Delta u|^2+\frac{|\nabla
u|^2}{|x|^2}+\frac{|u|^2}{|x|^4}\Big)dx\leq \epsilon.
\end{equation}
On the other hand, (\ref{equ44}) and Corollary \ref{co41} yields
\begin{equation*}
4\int|\Delta u|^2-4\int|u|^{2^{\#}}\geq
\tilde{C}_{\delta_0}\int|\Delta u_0|^2.
\end{equation*}
This and (\ref{equ85}) with $\epsilon=\epsilon_0\int|\Delta u_0|^2$,
implies that there exists $R_0>0$ such that for all $t\in [0,
\infty)$, we have
\begin{equation}\label{equ86}
4\int_{|x|\leq R_0}|\Delta u|^2-4\int_{|x|\leq R_0}|u|^{2^{\#}}\geq
C_{\delta_0}\int|\Delta u_0|^2.
\end{equation}
\begin{lemma}\label{lem81} Define \begin{equation}\label{equ87}
z_R(t)={\rm Im}\int
x\phi\big(\frac{|x|}{R}\big)\cdot\nabla\bar{u}udx, \end{equation}
where $\phi(r)$ is a smooth function with $\phi(r)=1$ when $r\leq1$
and $\phi(r)=0$ when $r\geq2$. Then
\begin{equation}\label{equ88}
z_R'(t)\geq4\int_{|x|\leq R}(|\Delta
u|^2-|u|^{2^{\#}})dx-O\Big(\int_{R\leq |x|\leq
2R}\big(\frac{|u||\Delta u|}{R^2}+\frac{|\nabla u||\Delta
u|}{R}+|\Delta u|^2\big)dx\Big).
\end{equation}
\end{lemma}
Thus, by (\ref{equ85}) and (\ref{equ86}), we get
\begin{equation*}
z'_R(t)\geq C_{d, \delta_0}\int|\Delta u_0|^2
\end{equation*}
for $R$ large.

Integrating in $t$, we have
\begin{equation*}
z_R(t)-z_R(0)\geq tC_{d, \delta_0}\int|\Delta u_0|^2.
\end{equation*}
But by the definition of $z_R(t)$, we have \begin{equation*}
|z_R(t)-z_R(0)|\leq 2R^4\|\Delta W\|_2^2,
\end{equation*}
which is a contradiction for $t$ large.

{\it Proof of Lemma \ref{lem81}.} We will compute
$\partial_tz_R(t)$.
\begin{align*}
\partial_tz_R(t)=&\ {\rm Im}\int
x\phi\big(\frac{|x|}{R}\big)\cdot\nabla\bar{u}_tudx+{\rm Im}\int
x\phi\big(\frac{|x|}{R}\big)\cdot\nabla\bar{u}u_tdx\\
=&\ 2{\rm Im}\int
x\phi\big(\frac{|x|}{R}\big)\cdot\nabla\bar{u}u_tdx+{\rm
Im}\int\nabla\cdot(x\phi\big(\frac{|x|}{R}\big))\bar{u}u_t\\
\doteq&\ A+B.
\end{align*}
$A$ and $B$ can be computed as follows:
\begin{align*}
A=&\ (4-d)\int\phi\big(\frac{|x|}{R}\big)\Big(|\Delta u|^2-|u|^{2^{\#}}\Big)dx\\
&\ +2{\rm Re}\int\phi'\big(\frac{|x|}{R}\big)\frac{(3-d)x\cdot\nabla
\bar{u}\Delta u}{R|x|}dx+2\int\phi'\big(\frac{|x|}{R}\big)|\Delta
u|^2dx\\
&\ +2{\rm Re}\int\phi''\big(\frac{|x|}{R}\big)\frac{x\cdot\nabla
\bar{u}\Delta
u}{R^2}dx+\frac{d-4}{d}\int\phi'\big(\frac{|x|}{R}\big)\frac{|x||u|^{2^{\#}}}{R}dx.
\end{align*}
\begin{align*}
B=&\ d\int\phi\big(\frac{|x|}{R}\big)\Big(|\Delta
u|^2-|u|^{2^{\#}}\Big)dx+2d{\rm
Re}\int\frac{x\cdot\nabla\bar{u}\Delta u}{R|x|}dx\\
&\ +d{\rm
Re}\int\Big(\phi''\big(\frac{|x|}{R}\big)\frac{\bar{u}}{R^2}+\phi'\big(\frac{|x|}{R}\big)\frac{(d-1)x\cdot\nabla\bar{u}}{R|x|}\Big)\Delta
udx+\int\phi'\big(\frac{|x|}{R}\big)\frac{|x||\Delta
u|^2}{R}dx\\
&\ +{\rm
Re}\int\Big(\phi'\big(\frac{|x|}{R}\big)\frac{d-1}{R|x|}+\phi''\big(\frac{|x|}{R}\big)\frac{d+1}{R^2}+\phi'''\big(\frac{|x|}{R}\big)\frac{|x|}{R^3}\Big)\bar{u}\Delta
udx\\
&\ +2{\rm
Re}\int\Big(\phi'\big(\frac{|x|}{R}\big)\frac{1}{R|x|}+\phi''\big(\frac{|x|}{R}\big)\frac{1}{R^2}\Big)x\cdot\nabla\bar{u}\Delta
udx.
\end{align*}
Since $\phi'$, $\phi''$ and $\phi'''$ are supported in $\{x:
R\leq|x|\leq 2R\}$, (\ref{equ88}) follows.

\vskip0.5cm
 \textbf{Acknowledgements:} C. Miao and G.Xu
were partly supported by the NSF of China (No.10725102,
No.10726053), and L.Zhao was supported by China postdoctoral science
foundation project.



\end{document}